\documentclass{amsart}
\usepackage{amsmath, amsthm, amscd, amsfonts}

\newcommand{\CC}{{\mathbb C}}

    \newtheorem{prop}{Proposition}[section]
    \newtheorem{cor}[prop]{Corollary}
    \newtheorem{lem}[prop]{Lemma}
    \newtheorem{rem}[prop]{Remark}
    \newtheorem{defn}[prop]{Definition}
    \newtheorem{thm}{Theorem}
    \newtheorem{example}{Example}
    \newtheorem{prob}{Problem}
  \renewcommand{\t}[1]{\mbox{$\tilde{#1}$}}
 \newcommand{\ip}[2]{\mbox{$(#1|#2)$}}
 
  \newcommand{\jbst}{$JB^*$-triple}
 \newcommand{\jbwst}{$JBW^*$-triple}
\newcommand{\jwst}{$JW^*$-triple}

 \newcommand{\csa}{$C^*$-algebra}
 
\newcommand{\jcsa}{$JC^*$-algebra}
 \newcommand{\jcst}{$JC^*$-triple}
  \newcommand{\pf}{{\bf Proof}.}
 
\newcommand{\ee}{\mbox{${\bf {\mathcal E}}$}}
  \newcommand{\sss}{\mbox{${\bf {\mathcal S}}$}}

   \newcommand{\tp}[3]{\{#1#2#3\}}
   \newcommand{\tpc}[3]{\{#1,#2,#3\}}
\begin{document}
\title[Contractive projections and operator spaces]
{Contractive projections and operator spaces}
\author{Matthew Neal \and Bernard Russo}
\thanks{This work was supported in part by NSF
 grant DMS-0101153}
\address{Department of Mathematics\\
        University of California, Irvine, California 92697-3875}
\email{mneal@math.uci.edu brusso@math.uci.edu}
\begin{abstract}
Parallel to the study of finite dimensional 
Banach spaces, there is a growing interest
in the corresponding local theory of operator spaces.
We define a family of Hilbertian operator spaces $H_n^k$,
$1\le k\le n$, generalizing the row and column Hilbert spaces $R_n,C_n$ and
show that an atomic subspace $X\subset B(H)$ which is
the range of a contractive projection on $B(H)$
is isometrically completely contractive
to an $\ell^\infty$-sum of the $H_n^k$ and 
Cartan
factors of types 1 to 4.
In particular, for finite dimensional $X$, 
this 
 answers  a question posed by Oikhberg and Rosenthal. Explicit in the proof
is a classification up to {\it complete isometry} of atomic w$^*$-closed
\jwst s without an infinite dimensional rank 1 w$^*$-closed ideal.
\end{abstract}
\maketitle


\begin{center}
{\bf Introduction}
\end{center}


It was shown by Choi-Effros  that an injective
operator system is  isometric to a conditionally complete 
\csa\ \cite[Theorem 3.1]{ChoEff77}. The fact that
 an injective operator system is the same as the image of a 
completely positive unital projection on $B(H)$
prompted a search for some algebraic structure in the range of 
a positive projection, or of a contractive projection.
A special case of a result of Effros-Stormer
showed that if a projection 
on a unital
\csa\ is positive and unital,
 then the range is isometric to
a Banach Jordan algebra \cite[Theorem 1.4]{EffSto79}.
Arazy-Friedman \cite{AraFri78} classified, up to Banach isometry, and
in Banach space terms,
 the range of an arbitrary contractive projection on the \csa\ of
all compact operators on a separable Hilbert space.
A special case of a result of Friedman-Russo showed  
that if a projection on a \csa\ is contractive, then the range 
is isometric to a Banach Jordan triple system
\cite[Theorem 2]{FriRus85}. 
Kaup \cite{Kaup84} extended the Friedman-Russo
result to contractive projections on \jbst s.

A consequence of these results is that, up to isometry, the ranges of the 
various projections can be
classified modulo a classification theorem of the various algebraic
structures involved. 
Recently, the operator space structure of the 
range of a {\it completely} contractive projection
has been studied. For projections acting on $B(H)$, such spaces coincide
with injectives in the category of operator spaces.
Christensen and Sinclair \cite[Theorem 1.1]{ChrSin89} prove
that every injective von Neumann algebra with separable predual which is not
finite type I of bounded degree is completely boundedly isomorphic to $B(H)$.
Robertson and Wasserman \cite[Corollary 7]{RobWas89} prove that an infinite
dimensional injective operator system on a separable Hilbert space is
completely boundedly isomorphic to either $B(H)$ or $\ell^\infty$.  
Robertson
and Youngson \cite[Theorem 1]{RobYou90} prove  that 
every injective operator space
is Banach isomorphic to one of $B(H),\ell^\infty,\ell^2$ or 
to a direct sum of these spaces.
Robertson
\cite[Corollary 3]{Robertson90} proves that an injective operator space
which is isometric to $\ell^2$ is completely isometric to $R$ or $C$
where $R$ and $C$ denote the row and column operator space versions of
$\ell^2$. These results can be thought of as giving a partial
classification
of injectives up to various types of isomorphisms. 

Note that the word injective in these examples is what we call 
1-injective below. Also, the spaces appearing in the above results are all
examples of
atomic \jwst s, but only 
$\ell^\infty$ and $B(H)$ are \csa s. Moreover, $B(H),R$ and $C$ are
examples of Cartan factors, while $\ell^\infty$ is a direct sum of 
countably many copies of the trivial Cartan factor $\CC$.

Operator spaces, that is, linear subspaces of $B(H)$, are the
appropriate setting for these types of problems. They were first studied
systematically in the thesis of Ruan \cite{Ruan88} and have been 
developed extensively since then
by Effros, Ruan, Blecher, Paulsen, Pisier, and others. 
Ruan \cite[Theorem 4.5]{Ruan89} showed that
an operator space is injective
if and only 
if it is completely isometric to $pAq$ for some injective \csa\ $A$ and
projections $p,q\in A$. Youngson had shown  earlier that
the range of a completely contractive projection on a $C^*$-algebra is
completely
isometric to a ternary algebra, that is, a subspace of a \csa\ that
is closed under the triple product $ab^*c$ \cite[Corollary~ 1]{Youngson83}.

Except for \cite{AraFri78},
there seem to be no results in the literature that 
classify the range of a contractive projection up to Banach isometry, or
up to completely bounded isomorphism. 
In this paper, we remedy this by investigating the 
structure of operator spaces which are 
the range of a contractive projection on
$B(H)$. These are known as 1-mixed injectives in operator space parlance.
We provide in Theorem~\ref{thm:2} a classification up to 
isometric complete contraction of 1-mixed injectives which are atomic.
In Theorem~\ref{thm:3}, we classify up to complete isometry all atomic
 w$^*$-closed \jwst s without an infinite dimensional rank 1 w$^*$-closed ideal. 
As a corollary, we show that an atomic (in particular, finite dimensional) contractively complemented subspace of
a $C^*$-algebra is a 1-mixed injective, that is, the range of a contractive projection on some $B(H)$.
Most
of these results have been announced in
\cite{NeaRus00}.

\section{Preliminaries}

An {\it operator space} is a subspace $X$ of $B(H)$, the space of
bounded linear operators on a complex Hilbert space. Its {\it 
operator
space structure} is given 
by the sequence of norms on the set of matrices $M_n(X)$
with entries from $X$, determined by the identification $M_n(X)
\subset M_n(B(H))=B(H\oplus H
\oplus \cdots \oplus H)$. For the basic theory of operator spaces 
and completely bounded maps,  we refer
to \cite{BlePau92},\cite{EffRua88},\cite{EffRua00}, \cite{Pisier1}, and \cite{Pisier2}
and the references therein. Let
us just recall that a linear 
mapping $\varphi:X\rightarrow Y$
between two operator spaces is {\it completely bounded}
if the induced mappings $\varphi_n:M_n(X)\rightarrow M_n(Y)$ defined by
$\varphi_n([x_{ij}])=[\varphi(x_{ij})]$
satisfy $\|\varphi\|_{\mbox{cb}}:=\sup_n\|\varphi_n\|<\infty$.
A completely bounded map is a {\it 
completely bounded isomorphism} if its inverse
exists and is completely bounded. Two operator spaces are 
{\it completely isometric}
if there is a linear isomorphism $T$ between them with $\|T\|_{\mbox{cb}}=
\|T^{-1}\|_{\mbox{cb}}=1$. We call $T$ a {\it complete isometry} in this case.

\medskip

In the matrix representation for $B(\ell^2)$ consider the {\it column
Hilbert space} $C=\overline{\mbox{sp}}\{e_{i1}:i\ge 1\}$ and the {\it row
Hilbert space} $R=\overline{\mbox{sp}}\{e_{1j}:j\ge 1\}$ and their finite 
dimensional versions
$C_n=\mbox{sp}\{e_{i1}:1\le i\le n\}$ and
$R_n=\mbox{sp}\{e_{1j}:1\le j\le n\}$. Here of course $e_{ij}$ is the 
operator defined by the 
matrix
with a 1 in the $(i,j)$-entry and zeros elsewhere. Although $R$ and $C$ are Banach isometric, they are not
completely isomorphic; and $R_n$ and $C_n$, while completely isomorphic,
 are not completely isometric.

An operator space $Z$ is {\it injective} if for any operator
space $Y$ and closed subspace $X\subset Y$, every completely bounded linear
map $T:X\rightarrow Z$ has a completely bounded extension $\t{T}:Y\rightarrow
Z$. In this case, there is a constant $\lambda\ge 1$
such that $\|\t{T}\|_{cb}\le \lambda\|T\|_{cb}$ and $Z$ is said to
be $\lambda$-injective.
If $\lambda=1$, then $Z$ is also called
 {\it isometrically injective}. 
A fundamental theorem in operator space theory is that $B(H)$
is 1-injective. This is the celebrated Arveson-Wittstock Hahn-Banach Theorem, see \cite[section 3]{EffRua88}.
It follows
that an operator space
$X\subset B(H)$ is $\lambda$-injective if and only if there is a completely
bounded projection $P$ from $B(H)$
onto $X$ with $\|P\|_{cb}\le \lambda$.

The literature on injective operator spaces cited in the introduction 
involves Cartan factors.
Cartan factors appeared 
 in the classification of   Jordan triple systems
and bounded symmetric domains. There are six types of Cartan factors of which
four will be relevant to our work. 
A Cartan factor of type 1 is
the space $B(H,K)$ of all bounded operators from one complex Hilbert space
$H$ to another $K$.  By fixing orthonormal bases for $H$ and $K$, we may
think of $B(H,K)$ as all $\dim K$ by $\dim H$ matrices which define
bounded operators.
To define the Cartan factors of types 2 and 3 we need to fix
a conjugation $J$ on a Hilbert space $H$, that 
is, a conjugate-linear isometry of order 2. 
Then a Cartan factor of type 2 (respectively type 3) is $A(H,J)=\{
x\in B(H):x^t=-x\}$ (respectively $S(H,J)=\{x\in B(H):x^t=x\}$), where
$x^t=Jx^*J$. Since conjugations
are in one-to-one correspondence with orthonormal bases of $H$, 
we may think of these as anti-symmetric (resp. symmetric)
$\dim H$ by $\dim H$ matrices which define bounded operators.
A Cartan factor of type 4, or {\it spin factor} 
 will be described in more detail
in subsection~\ref{ss:spin}.  

The following concepts were 
 introduced by Oikhberg and Rosenthal 
in \cite[section 3]{OikRospp} in their study of 
extension properties for the space of compact operators.
The operator space $Z$ is a {\it mixed injective} if for every completely
bounded linear map $T$ from an operator space
$X$ into $Z$ and any operator space $Y$ containing $X$, $T$
 has a bounded extension $\tilde{T}$
to $Y$. 
In this case, there is a constant $\lambda\ge 1$
such that $\|\t{T}\|\le \lambda\|T\|_{cb}$,
 and $Z$ is said to
be $\lambda$-mixed injective.
A
$1$-mixed injective operator space is also said to be {\it isometrically
mixed injective} and  $X$ is $\lambda$-mixed injective if and only if there is a bounded
projection $P$ from $B(H)$
onto $X$ with $\|P\|\le \lambda$.
An operator space $X$ is {\it completely semi-isomorphic} to an operator
space $Y$ if there is a linear homeomorphism $T:X\rightarrow Y$ which
is completely bounded. Such a $T$ is called a {\it complete semi-isomorphism}.
If in addition $\|T\|_{cb}=\|T^{-1}\|=1$, then
$X$ is {\it completely semi-isometric} to $Y$ and $T$ is a {\it
complete semi-isometry}. It is shown in 
\cite[Proposition 3.9]{OikRospp}
that mixed injectivity is preserved by complete semi-isomorphisms in the sense
that if $Y$ is a mixed injective, then so is $X$.

The Cartan factors of types 1 to 4 are examples of 1-mixed injectives. This
is obvious for types 1,2,3 and for type 4 it is proved in \cite[Lemma 2.3]{EffSto79}.
Cartan factors of types 5 and 6 will play no role in this paper since 
neither is
 even isometric to a 1-mixed injective operator space.
For if it were, it would follow from  \cite[Theorem 2]{FriRus85}
that it would
 be isometric to a \jcst\ (defined below). This is impossible since they
are well-known to be ``exceptional'' 
( i.e., not triple isomorphic to the Jordan triple structure induced by
an associative $^*$-algebra, see \cite[2.8.5]{HanSto84} for the Jordan
algebra version of this),  
and surjective isometries coincide with
triple isomorphisms (the latter is proved for  \jcst s in \cite{Harris73}; for the more general class of \jbst s see
\cite{Kaup83} or \cite[Lemma 1]{Bar1}). The space of compact operators on a separable Hilbert space
is not a 1-mixed injective, but it seems to be an interesting open question
whether it has the {\it mixed separable extension property} \cite{OikRospp}, 
that is, in the 
definition of mixed injective, only separable operator spaces $X\subset Y$
are considered.

In view of the relaxed definition of 1-mixed injectives, one cannot 
immediately expect a classification of them up to complete
isometry (however, see Theorem~\ref{thm:3}). It is more natural to ask 
for a classification of 1-mixed injectives up to complete 
semi-isometry.  In order to
formulate our results precisely we recall some basic facts about 
\jcst s.

A \jcst\ 
is a norm 
closed complex linear subspace $M$ of a $C^*$-algebra $A$
which is closed under the operation $a\mapsto aa^*a$. 
\jcst s were defined and studied (using the name
$J^*$-algebra)
 as a generalization of \csa s by Harris \cite{Harris73} in
connection with function theory on infinite dimensional bounded symmetric
domains. By a polarization identity, any \jcst\ is closed under the
triple product
\begin{equation}\label{eq:product}
(a,b,c)\mapsto
\tp{a}{b}{c}:=\frac{1}{2}(ab^*c+cb^*a),
\end{equation}
under which it becomes a Jordan triple system.  In this paper, the notation
$\tp{a}{b}{c}$ will always denote the triple product (\ref{eq:product}). 
A linear map which preserves the triple product (\ref{eq:product}) will be
called a {\it triple homomorphism}.
Cartan factors are 
examples of \jcst s, as are \csa s, and Jordan \csa s.

A {\it JW*-triple} is defined to be a \jcst\ which is a dual space. It follows from
\cite[Corollary 9]{Bar1} that a \jwst\ is isometric to a \jcst\ which is weak operator closed.

Note that some of the results about \jcst s that we are going to cite were
proven for the more general class of \jbst s. For example, 
\cite[Theorem 2.1]{Bar2} shows that all preduals of a \jwst\ are isometric. 
\jbst s, in and of themselves,
 will play no role in this paper, but the interested reader can
consult \cite{Russo94} for a comprehensive 
survey from an operator algebra point
of view.

A special case of a \jcst\ is a {\it
ternary algebra}, that is, a subspace of $B(H,K)$ closed under the
{\it ternary product} $(a,b,c)\mapsto ab^*c$. A {\it ternary homomorphism}
is a linear map $\phi$ satisfying $\phi(ab^*c)=\phi(a)\phi(b)^*\phi(c)$. 
These spaces are also 
called, more appropriately, 
{\it associative triple systems}.
They have been studied both concretely in \cite{Hestenes62}
and abstractly in \cite{Zettl83}. We shall use the term ternary 
algebra in this paper, but we shall not need any special results about them, other
than the well-known and simple fact that a ternary isomorphism between two ternary
algebras is a complete isometry.
A key step in our proof of Theorem~\ref{thm:2}
 will be to extend
a Banach 
isometry between two \jwst s to a ternary isomorphism of
their ternary envelopes.

If $v$ is a partial isometry in a \jcst\  $M\subset B(H,K)$, 
then the projections
$l=vv^*\in B(K)$ and $r=v^*v\in B(H)$ give rise to (Peirce)
projections $P_k(v):M\rightarrow M,\ k=2,1,0$  as
follows; for $x\in M$,
\[
P_2(v)x=lxr\quad , \quad P_1(v)x=lx(1-r)+(1-l)xr\quad , \quad
P_0(v)x=(1-l)x(1-r).
\]
These projections $P_k(v)$
are easily seen to have the following properties. They
are contractive projections and their ranges, denoted by $M_k(v)$ are
$JC^*$-subtriples of $M$ satisfying 
$M=M_2(v)\oplus M_1(v)\oplus M_0(v)$. They obey
{\it Peirce calculus}, by which is meant
\[
\tp{M_2(v)}{M_0(v)}{M}=\tp{M_0(v)}{M_2(v)}{M}=0\quad , \quad
\tp{M_i(v)}{M_j(v)}{M_k(v)}\subset M_{i-j+k}(v)
\]
where it is understood that $M_{i-j+k}(v)=\{0\}$ if
$i-j+k\not\in \{0,1,2\}$.

The Peirce space $M_2(v)$ plays a special role. It has the structure of
a unital Jordan $^*$-algebra with unit $v$
under the product $(a,b)\mapsto a\circ b:=\tp{a}{v}{b}$
and involution $a\mapsto a^\sharp:=\tp{v}{a}{v}$. 
For example, the Jordan identity
$(a\circ a)\circ (a\circ b)=a\circ ((a\circ a)\circ b)$ amounts to 
$\tp{a}{v}{\tp{\tp{a}{v}{a}}{v}{b}}=\tp{\tp{a}{v}{a}}{v}{\tp{a}{v}{b}}$,
which is trivial to verify for \jcst s. For more general Jordan triple
systems see for example 
\cite[3.13]{Loos77},\cite{Loos75}, or 
\cite[19.7]{Upmeier85},
which are references for the general
theory of (Banach) Jordan triple systems. 

We shall write $M_2(v)^{(v)}$ to
denote the space $M_2(v)$ with this structure. If $M=M_{2}(v)$, then we refer to $M_2(v)^{(v)}$ as 
an {\it isotope} 
 of $M$. If $M$ is
a ternary algebra, then $M_2(v)^{(v)}$ is a unital \csa\ with product 
$a \cdot b=av^{\ast}b$, involution $a^{\sharp}=va^{\ast}v$, and
unit $v$. In this case, the identity map
from  $M_{2}(v)$ to $M_{2}(v)^{(v)}$ is a ternary isomorphism, since $ab^*c
=av^*b^\sharp v^*c$, and hence also a complete isometry.


\medskip  

A partial isometry $v$ is said to be {\it 
minimal} in $M$ if $M_2(v)=\CC v$. This
is equivalent to $v$ not being the sum of two orthogonal non-zero 
partial isometries.
Recall that two partial isometries $v$ and $w$
(or any two Hilbert space operators) are
orthogonal if $v^*w=vw^*=0$. This is equivalent to $v\in M_0(w)$ and
will be denoted by $v\perp w$.
Each finite dimensional \jcst\ is
the linear span of its minimal paritial isometries. More generally, an
{\it atomic} \jwst\ is defined to be one which is the weak$^*$-closure of
the span of its minimal partial isometries.  The {\it rank} of a \jcst\ is the maximum
number of mutually orthogonal minimal partial isometries. For example, the
rank of the Cartan factor of type 1 $B(H,K)$ is the minimum of the dimensions
of $H$ and $K$; and the rank of the Cartan factor of type 4 is 2. Other
relations between two partial isometries that we shall need are defined 
in terms of the Peirce spaces as follows.  Two partial isometries $v$ and
$w$ are said to be {\it colinear} if $v\in M_1(w)$ and $w\in M_1(v)$, notation
$v\top w$.
A partial isometry $w$ is said to {\it govern} $v$ if $v\in M_2(w)$ and
$w\in M_1(v)$.  It is easy to check that $v\in M_j(w)$ if and only if
$\tp{w}{w}{v}=(j/2)v$, for $j=0,1,2$. 

\jcst s of arbitrary dimension occur naturally in functional 
analysis and in holomorphy. 
As noted in the introduction, a
 special case of a theorem of Friedman and Russo \cite[Theorem 2]{FriRus85}
 states that if $P$ is a contractive
projection on a $C^*$-algebra $A$, then there  is a linear isometry
of the range $P(A)$ of $P$ onto a $JC^*$-subtriple of $A^{**}$. 
A special case of a theorem of Kaup \cite{Kaup83} gives a bijective correspondence between Cartan factors and
irreducible bounded symmetric domains in complex Banach spaces.

There is a structure theorem for
atomic \jwst s, for which we refer to \cite[p.\ 302]{DanFri87},\cite{Horn},
\cite{Neher86} for the proofs. A \jwst\ is {\it irreducible} if it is not the $\ell^\infty$-direct sum of 2 non-zero
w$^*$-closed ideals. 
The version of the structure theorem that  we shall use is the following.

\begin{lem} \label{atomic}
 Each atomic \jwst\ $X$ is the 
$\ell^\infty$-direct sum $X=\oplus^{\ell^\infty}_\lambda X_\lambda$ of weak-$^*$ closed irreducible ideals, 
and each summand $X_\lambda$ is the weak$^*$-closure 
of the complex linear span of a {\it grid} of minimal partial
isometries. 
Grids come in four types and each $X_\lambda$ is Banach
isometric and hence triple isomorphic
 to a Cartan factor of one of the types 1--4. 
\end{lem}

We shall describe the grids
for the Cartan factors of types 1--4 (the so-called rectangular grid, 
symplectic grid, hermitian grid, 
and spin grid) when they are needed later in this
paper.  As will be seen,
grids only give information about
the symmetrized triple product (\ref{eq:product}), whereas the operator
space structure depends on the ternary product 
$(a,b,c)\mapsto ab^*c$. 

It follows from Lemma~\ref{atomic} and \cite[Theorem 2]{FriRus85}
 that a finite dimensional 1-mixed injective operator
space is Banach isometric to an $\ell^\infty$-direct sum of Cartan factors.
Oikhberg and Rosenthal \cite[Problem 3.3]{OikRospp} ask 
whether every finite dimensional
1-mixed injective operator space is in fact completely semi-isometric to an
$\ell^\infty$-direct sum of Cartan factors of types 1-4. 
Corollary~\ref{cor:2.1} of Theorem~\ref{thm:2}
below
answers  this question.

 Theorem~\ref{thm:2} below is formulated for {\it atomic} 1-mixed injective
operator spaces. A Banach space $X$ with predual $X_{\ast}$ is said to be
atomic if the closed
unit ball $X_{\ast,1}$ is the norm closed convex hull of its extreme points. In particular, reflexive
Banach spaces are atomic as are the duals of unital \csa s. 

\section{Main Results and Reduction}

In this section we state Theorems~\ref{thm:1},\ref{thm:2}, and \ref{thm:3},
and give a reduction for the proof of Theorem~\ref{thm:2}.

\begin{thm}\label{thm:1}
There is a family of 
1-mixed injective Hilbertian operator spaces 
$H_n^k$, $1\le k\le n$, of finite dimension $n$,
with the following properties:

\begin{description}
\item[(a)] $H_n^k$ is a  subtriple  of the Cartan factor of type 1 
consisting of
all $n\choose k$ by $n\choose n-k+1$ complex matrices.
\item[(b)] Let $Y$ be a \jwst\ of rank 1 {\rm (}necessarily 
atomic {\rm )}.
\begin{description}
\item[(i)] If $Y$ is of finite dimension $n$ then it is isometrically
completely contractive to some
$H_n^k$.  
\item[(ii)] If $Y$ is infinite dimensional then it is isometrically completely
contractive to $B(H,\CC)$ or $B(\CC,K)$.
\end{description} 
\item[(c)] $H_n^n$ {\rm (}resp. $H_n^1${\rm )}
coincides with $R_n$ {\rm (}resp. $C_n${\rm )}.
\item[(d)] For $1<k<n$, $H_n^k$ is not completely semi-isometric to $R_n$
or $C_n$.
\end{description}
\end{thm}

The spaces $H_n^k$ are explicitly constructed in section 6. These spaces 
appeared
in a slightly different form in \cite{AraFri78}, see Remark~\ref{rem:6.17}.
The authors are grateful to N. Ozawa for showing us the proof of (d).

\begin{thm}\label{thm:2}
Let $X$ be a 1-mixed injective operator space which is atomic.
Then $X$ is completely semi-isometric 
to a direct sum of Cartan factors of types 1 to
4 and the spaces $H_n^k$.
\end{thm}

The following Corollary to Theorem~\ref{thm:2}, together with (d) of
Theorem~\ref{thm:1},
answers  the question of Oikhberg and Rosenthal
\cite[Problem 3.3]{OikRospp}. 

\begin{cor}\label{cor:2.1}
A finite dimensional 1-mixed injective operator space is
 completely semi-isometric 
to a direct sum of Cartan factors of types 1 to
4 and the spaces $H_n^k$.
\end{cor}

We now begin the proofs of Theorems~\ref{thm:1} and ~\ref{thm:2}. 

Let $X\subset B(H)$ be a 
1-mixed injective
operator space. Then there is a contractive projection
on $B(H)$ with range $X$. 
By \cite[Theorem 2]{FriRus85}, there is thus
a linear isometry $\ee_0$ from $X$ onto a \jcst\ $Y\subset
A:=B(H)^{**}$ of the form
$\ee_0(x)=pxq$ for suitable projections $p,q$ in the von Neumann algebra
$A$. Since $(\ee_0)_n:M_n(X)\rightarrow M_n(Y)$ has the form 
$[x_{ij}]\mapsto \mbox{diag}\,  (p,p,\ldots,p)\ [x_{ij}]\ \mbox{diag}\, (q,q,\ldots,q)$
we have  the following lemma.

\begin{lem}\label{lem:3.1}
$\ee_0$ is completely contractive and hence a complete semi-isometry of
$X$ onto the \jcst\ $Y$.
\end{lem}

\begin{lem}\label{batomic}
Suppose $X$ is a Banach space with predual which is isometric to a 
$JC^*$-triple $Y$. Then  $X$ is atomic as a Banach space
if and only if $Y$ is an atomic JW*-triple.
\end{lem}
\pf\
Since $X_{\ast}$ is a predual of $Y$, $Y$ is a JW*-triple. Assume $X$ is atomic as a Banach space. It is shown in \cite[Prop
4c]{FriRus85bis} that the minimal partial isometries
$v$ in $Y$ are in 1-1 correspondence with extreme points $\phi$ of $Y_{\ast,1}$ via the mapping $\phi \rightarrow
v$ if
$\phi(v)=1$. By \cite[Theorem 2]{FriRus85bis} $Y$ has an internal $\ell^\infty$ direct sum decomposition ${\mathcal A}
\oplus {\mathcal N}$ into w*-closed
subtriples where
$\mathcal A$ is and atomic and $\mathcal N$ contains no minimal partial isometries. It follows that ${\mathcal N} = \{ 0 \}$ and $Y$ is
atomic. 

Conversely, if $Y$ is an atomic \jwst, then by \cite[Theorem 1]{FriRus85bis}, $Y_*=A\oplus^{\ell^1} N$, where $N$
has no extreme points.  It follows that $N=\{0\}$.
\qed

\begin{lem}\label{lem:3.2}
It suffices to prove Theorem~\ref{thm:2} in the case that $Y$ 
{\rm (}$=\ee_0(X)${\rm )} is 
triple isomorphic to a Cartan factor.
\end{lem}
\pf\
By Lemma \ref{atomic} and Lemma \ref{batomic}, $Y=\oplus_\alpha
Y_\alpha$ is the internal $\ell^\infty$
direct sum of a family of subtriples $Y_\alpha$, 
each of which is triple isomorphic
to a Cartan factor  of one of the types 1--4.

Suppose that $T_\alpha:Y_\alpha\rightarrow Z_\alpha$ is a complete semi-isometry.
Then $\oplus T_\alpha:\oplus Y_\alpha\rightarrow \oplus Z_\alpha$ is also 
a complete semi-isometry, by the following commutative diagram:
\[
\begin{array}{ccccc}
M_n(\oplus Y_\alpha)&\stackrel{\mbox{isometry}}{\longrightarrow} & \oplus M_n(Y_\alpha)\\
                    &&\\
(\oplus T_\alpha)_n\downarrow\quad\quad\quad & &\quad\quad\quad\downarrow\oplus (T_\alpha)_n\\
&&\\
M_n(\oplus Z_\alpha)&\stackrel{\mbox{isometry}}{\longrightarrow} & \oplus M_n(Z_\alpha)\\
\end{array}.
\]

To show the isometry part of the above diagram, one can use the idea of 
\cite[Lemma 1.3]{FriRus85bis}. 
For completeness, we include the argument. For $c\in M_n(X\oplus Y)$ with $c_{ij}=a_{ij}\oplus b_{ij}$, we have
$c=[c_{ij}]=[a_{ij}\oplus 0]+[0\oplus b_{ij}]=a+b$ with $a,b$ orthogonal operators, that is, $ab^*=a^*b=0$. Then, assuming
$\|a\|\le 1$ and $\|b\|\le 1$,
$\|c\|=\|a+b\|=\|(a+b)(a+b)^*(a+b)\|^{1/3}=\|a^{3^n}\|^{3^{-n}}=(\|a^{3^n}+b^{3^n}\|^{3^{-n}})^{3^{-n}}\le
\|a\|^{3^n}+\|b\|^{3^n}\le 2^{3^{-n}}\rightarrow 1$, proving that $\|c\|_{M_n(X\oplus Y)}\le \|[a_{ij}] \oplus
[b_{ij}]\|_{M_n(X)\oplus M_n(Y)}$.  Conversely, assume $\|a\|=1$.  Then $1=\|a\|^5=\|aa^*aa^*a\|=\|aa^*(a+b)a^*a\|\le \|a+b\|$,
so that  
$\|a\|\le \|c\|$ and it follows that $\|[a_{ij}] \oplus [b_{ij}]\|_{M_n(X)\oplus M_n(Y)} \leq \|c\|_{M_n(X\oplus Y)}$.
\qed

\medskip

An {\it ideal} of a $JC^*$-triple $Y$ is a subspace $I\subseteq Y$ 
such that
$\{ Y \,\ I \,\ Y \} +  \{ I \,\ Y \,\ Y \} \subseteq I$.
By \cite[Prop.\ 2.1]{FriRus83}, 
 the second dual of a $JC^*$-triple is a $JW^*$-triple. 

\begin{lem}\label{p}
Every $JW^*$-subtriple $Y$ of a $C^*$-algebra $A$ is 
completely semi-isometric to a w*-closed $JW^*$-subtriple of
$A^{\ast\ast}$.
\end{lem}
\pf\
 By separate 
w*-continuity of multiplication, the annihilator
$Y_{\ast}^{0}$ is  a w*-closed ideal of $Y^{\ast\ast}$. By
 \cite{Horn} or \cite[Theorem 3.5]{Neher86}, $Y^{**}=Y_*^0\oplus^{\ell^\infty} J$,
where $J$ is a w*-closed ideal orthogonal to $Y_*^0$.
Let $P$ (resp. $Q$) be the projection of $Y^{**}$ onto $Y_*^0$
(resp. $J$). For each element $z\in Y^{**}$, there is $y\in Y$ with
$z-y\in Y_*^0$.  It follows that $Q(Y)=Q(Y^{**})=J$, and it is easy
to see by the orthogonality that $Q$ is a w*-continuous triple homomorphism
from $Y^{**}$ onto $J$.  Since $P(Y^{**})=Y_*^0$, $Q$ is one-to-one on
$Y$. As in the proof of Lemma~\ref{lem:3.2}, $Q=0\oplus Id_{J}$ 
is a complete contraction.
\qed 

\medskip

Theorems~\ref{thm:1}(b) and ~\ref{thm:2} are immediate
consequences
 of Lemmas \ref{lem:3.1}, \ref{lem:3.2}, \ref{p} and the following proposition.
Theorem~\ref{thm:1}(a) is proved in section 6 (see Remark~\ref{rem:6.05}).
  Theorem~\ref{thm:1}(c) is proved in
Proposition~\ref{prop:5.10bis} and Theorem~\ref{thm:1}(d) is proved at 
the end of section 7.

\begin{prop}\label{prop:1} Let $Y$ be a \jwst\
which is w$^*$-closed in a $W^*$-algebra. If $Y$
is either of rank at least 2, or of rank 1 and infinite
dimensional, and is triple isomorphic
to a Cartan factor of type 1,2,3, or 4, then it is in fact completely
semi-isometric to a Cartan factor of the same type. A 
finite dimensional 
\jwst\ which is triple isomorphic to a Cartan factor of rank
1, is completely semi-isometric to one of the
spaces $H_n^k$.
\end{prop}

As a by-product of the proof of Proposition~\ref{prop:1}, we 
shall obtain the following theorem, which by Lemma~\ref{atomic} gives a classification 
up to {\it complete isometry}
of atomic \jwst s which are w$^*$-closed and
contain no infinite dimensional w$^*$-closed ideals of rank 1.

We need the following definitions. If $B(H,K)$ is a Cartan factor of
type 1, then 
\[
\mbox{Diag}\, (B(H,K),B(K,H)):=\{(x,x^t): x\in B(H,K)\},
\]
where the transpose is respect to fixed orthonormal bases for $H$
and $K$.  We give $\mbox{Diag}\, (B(H,K),B(K,H))$ the operator space
structure induced by its natural embedding in $B(H\oplus K\oplus H\oplus K)$ 
and note that $\mbox{Diag}\, (B(H,K),B(K,H))$ is contractively complemented 
therein. Indeed, first project $B(H\oplus K\oplus H\oplus K)$ onto
$B(H,K)\oplus B(K,H)$ and follow by $(x,y)\mapsto ((x+y^t)/2,(x^t+y)/2)$.

For a fixed dimension $n$ and 
each $j=1,\ldots,m$ let $H_j$ be a Hilbert space of 
dimension $n$ with a specified
orthonormal basis ${\mathcal B}_j=\{e_{j,1},\ldots,
e_{j,n}\}$. Then 
\[
\mbox{Diag}\, (\{H_j,{\mathcal B}_j\}):=\{(\sum_k\alpha_ke_{1k},\sum_k\alpha_ke_{2k},
\ldots,
\sum_k\alpha_ke_{mn}):\alpha_k\in\CC,1\le k\le m\}
\] 
The space $\mbox{Diag}\, (\{H_j,{\mathcal B}_j\})$ is contractively complemented in $\oplus_{j=1}^n H_j$.

\begin{thm} \label{thm:3}  Let $Y$ be an 
atomic w$^*$-closed $JW^*$-subtriple of a $W^*$-algebra.
\begin{description}
\item[(a)] If $Y$ is irreducible and
 of rank at least 2, then it  is completely isometric
to a Cartan factor of type 1--4 or the
space $\mbox{Diag}\, (B(H,K),B(K,H))$.
\item[(b)] If $Y$ is of finite dimension $n$
 and of rank 1, then it is completely
isometric to $\mbox{Diag}\, (H_n^{k_1},\ldots,H_n^{k_m})$, for appropriately
chosen bases defined in section 7, and where
$k_1>k_2>\cdots>k_m$.
\item[(c)] $Y$ is completely semi-isometric to a direct sum of the spaces in (a) and (b).
If $Y$ has no infinite dimensional rank 1 summand, then it is completely
isometric to a direct sum of the spaces in (a) and (b).
\end{description}
\end{thm}

\begin{cor}
Every finite dimensional \jcst\ is completely isometric to an $\ell^\infty$-direct sum of
Cartan factors of types 1--4 and the spaces $\mbox{Diag}\, (B(H,K),B(K,H))$ and
$\mbox{Diag}\, (H_n^{k_1},\ldots,H_n^{k_m})$.
\end{cor}

\begin{cor}
Every atomic contractively complemented subspace $X$ of a $C^*$-algebra $A$ is 1-mixed injective.
\end{cor}
\pf\
By \cite[Theorem 2]{FriRus85} and Lemma~\ref{batomic}, the map $\ee_0$ mentioned earlier
is a complete semi-isometry of $X$ onto an atomic $JW^*$-subtriple of $A^{**}$. By Theorem~\ref{thm:3}
and Lemma~\ref{p}, $X$ is completely semi-isometric to a direct sum of the spaces listed above, which
as noted are 1-mixed injectives.  Then by \cite[(3.9)]{OikRospp}, $X$ is 1-mixed injective.\qed

\begin{cor}
Every atomic JW*-triple is a 1-mixed injective.
\end{cor}

We shall prove 
Proposition~\ref{prop:1} and Theorem~\ref{thm:3}
case by case in the following sections. Cartan factors of Types 
3 and 4 are handled in
section 3, type 2 in section 4 and type 1 in sections 5,6 and 7.
Section 6 also introduces the spaces $H_n^k$ and section 7 also gives some examples and states
some open problems.

\section{Cartan factors of types 3 and 4}

The Cartan factors of types 3 and 4 have a unital Jordan $^*$-algebra
structure in which we frame the proofs of Proposition~\ref{prop:1} and
Theorem~\ref{thm:3}.

\subsection{Cartan factors of type 4}\label{ss:spin}

We first prove Proposition~\ref{prop:1} and Theorem~\ref{thm:3}
in the case that 
$Y$ is triple isomorphic to a Cartan factor of type 4. Let us first
describe the concrete model which we use for such a Cartan factor, from
\cite[Theorem 6.2.2]{HanSto84} and \cite{Harris73}.

A {\it spin system}  is a subset $\sss=\{1,s_1,\ldots,s_k\}$
of self-adjoint elements of $B(H)$ containing the unit  and 
satisfying $s_{i}s_{j}+s_{j}s_{i}=\delta_{ij}2$.
It follows that  
$\mbox{sp}_{\CC}\, \sss$ is a $(k+1)$-dimensional Jordan $C^*$-subalgebra of
$B(H)$. A {\it 
spin factor} is a subspace $X$ of $B(H)$
of dimension at least 2  which is the closed linear span of a spin system
of {\it arbitrary} cardinality. 

We now recall the standard matrix representation of the spin factor 
$\mbox{Sp}\, (n)$,
$3\le n\le \infty$, for the separable case (c.f. \cite[6.2.1]{HanSto84}),
which is the Cartan factor of type~ 4.
Let
\[
\sigma_1=\left[ \begin{array}{rr}
1&0\\ 0&-1
\end{array}\right]\quad , \quad
\sigma_2=\left[ \begin{array}{rr}
0&1\\ 1&0
\end{array}\right]\quad , \quad
\sigma_3=\left[ \begin{array}{rr}
0&i\\ -i&0
\end{array}\right]
\]
be the Pauli spin matrices. Denote by $\sigma_3^n$ the
$n$-fold tensor product $\sigma_3\otimes\cdots\otimes \sigma_3$ of $\sigma_3$
with itself $n$ times in $M_{2^n}(\CC)$.
Define
\[
s_1=\sigma_1\ ,\ s_2=\sigma_2\ ,\ s_3=\sigma_3\otimes\sigma_1\ ,\ s_4=
\sigma_3\otimes\sigma_2\ldots
\]
and in general $s_{2n+1}=\sigma_3^n\otimes \sigma_1$ and
$s_{2n+2}=\sigma_3^n\otimes \sigma_2$.

With the imbeddings $M_{2^n}(\CC)\subset M_{2^{n+1}}(\CC)$ given by
\[
a\mapsto a\otimes 1=\left[ \begin{array}{rr}
a&0\\ 0&a
\end{array}\right]
\]
we have $s_k\in M_{2^n}(\CC)$ if $k\le 2n$ and $\{1,s_1,\ldots,s_k\}$ is
a spin system for each $k\ge 2$.  The linear span $\mbox{Sp}\, (k+1)$ of 
$\{1,s_1,\ldots,s_k\}$ is a $(k+1)$-dimensional spin factor contained
in $M_{2^n}(\CC)$ if $2\le k\le 2n$. For more details and the case
$k=\aleph_0$, see \cite[Theorem 6.2.2]{HanSto84}.

As an operator space, a spin factor
$X$ is determined up to complete isometry by 
the cardinality  of the spin system. Indeed, it is easy to see that two finite
 spin systems with the same 
number of elements generate $C^*$-algebras which 
are $^*$-isomorphic with basis
consisting of all finite products of elements in the spin system 
and the unit. If a spin system $\{s_{\lambda}\}_{\lambda
\in \Lambda}$ has arbitrary cardinality, the inductive limit 
$\overline{\cup_F A(F)}$ (norm closure)
of the collection of $C^*$-algebras $A(F)$ generated by finite subsets $F$ of 
$\{s_{\lambda}\}_{\lambda \in \Lambda}$ is exactly the $C^*$-algebra generated by the spin system. Thus, any
two spin systems with the same cardinality generate *-isomorphic 
$C^*$-algebras
(cf. \cite[1.23]{Sakai} or \cite[Theorem 5.2.5]{BraRob81}). 

Suppose now that the $Y$ in the
statement of Proposition~\ref{prop:1}
 is triple isomorphic to a Cartan factor of type 4. Let $A$ denote
any von Neumann algebra containing $Y$.  The \jwst\ $Y$ contains 
a {\it spin grid} 
$\{u_j,\t{u}_j:j
\in J\}$, or
$\{u_j,\t{u}_j:j \in J\}\cup\{u_0\}$ in
the case that $Y$ is of finite odd dimension.

Let us recall the properties of a spin grid from \cite[p.\ 313]{DanFri87}. The elements
$u_j$ and $\t{u}_j$ 
(but not $u_0$) are minimal non-zero partial isometries; for $i\ne j$,
$u_i$ is colinear
with $u_j$ and with $\t{u}_j$, and $\t{u}_j$ is colinear with $\t{u}_i$; and
for $i\ne j$,
\begin{equation}\label{eq:quad}
\tp{u_i}{u_j}{\t{u}_i}=-\frac{1}{2}\t{u}_j\ ,\ 
\tp{u_j}{\t{u}_i}{\t{u}_j}=-\frac{1}{2}\t{u}_i.
\end{equation}
In case $u_0$ is present, for each $i\ne 0$, $u_0$ governs $u_i$ and
$\t{u}_i$, and
\begin{equation}\label{eq:govern}
\tp{u_0}{u_i}{u_0}=-\t{u}_i\ ,\ \tp{u_0}{\t{u}_i}{u_0}=-u_i.
\end{equation}
All other triple products from the spin grid are 0, and in particular,
$u_i$ is orthogonal to $\t{u}_i$.

It is not hard to see (c.f. \cite{DanFri87}) that the complex span $Y$ of a spin grid has an equivalent Hilbertian norm and is
hence reflexive. It is also clear from the grid properties that all such $Y$ are rank 2.

Let
$v=i(u_1+\t{u}_1)$, where 1 is an arbitrary element of the index set
$J$. It is easy to see that that $Y=Y_2(v)$. As noted in 
the preliminaries, $A_{2}(v)$ and $A_{2}(v)^{(v)}$ are
ternary isomorphic and thus completely isometric.
Thus, the identity map
$Y\rightarrow Y_2(v)^{(v)}$ is a complete isometry. 

The following lemma is easily verified by using (\ref{eq:quad}),
(\ref{eq:govern}) and Peirce calculus. For the convenience of the reader, 
we include some of the details.

\begin{lem}
$Y_2(v)^{(v)}$ is a Cartan factor of type 4. More precisely, let
$s_{j}=u_j+\t{u}_j,\ j \in J-\{1\};
\  t_{j}=i(u_j-\t{u}_j),\ j \in J$.
Then a spin system in the $C^*$-algebra 
$A_{2}(v)^{(v)}$ which linearly spans $Y_2(v)^{(v)}$ is given by
\[
\{s_j,t_k,v:j\in J-\{1\},k\in J\};
\]
or, if the spin factor is of odd finite dimension,
\[
\{s_j,t_k,v,u_0:j\in J-\{1\},k\in J\}
\]
\end{lem}
\pf\
If $j\ne 1,k\ne 1$,
\begin{eqnarray*}
s_j\cdot s_k+s_k\cdot s_j&=&s_jv^*s_k+s_kv^*s_j=2\tp{s_j}{v}{s_k}\\
&=&-2i\tpc{u_j+\t{u}_j}{u_1+\t{u}_1}{u_k+\t{u}_k}.
\end{eqnarray*}
If $j\ne k$ then all 8 terms in the expansion of this triple product are
zero since the triple product of three mutually colinear partial isometries
is zero.  On the other hand,
\begin{eqnarray*}
2s_j\cdot s_j&=&-2i\tpc{u_j+\t{u}_j}{u_1+\t{u}_1}{u_j+\t{u}_j}\\
&=&-2i[\tp{u_j}{u_1}{u_j}+\tp{u_j}{u_1}{\t{u}_j}+\tp{u_j}{\t{u}_1}{u_j}+
\tp{u_j}{\t{u}_1}{\t{u}_j}]\\
&+&-2i[\tp{\t{u}_j}{u_1}{u_j}+\tp{\t{u}_j}{u_1}{\t{u}_j}+
\tp{\t{u}_j}{\t{u}_1}{u_j}+\tp{\t{u}_j}{\t{u}_1}{\t{u}_j}]\\
&=&-2i[0-\t{u}_1/2+0-u_1/2-\t{u}_1/2+0-u_1/2+0]\mbox{ (by (\ref{eq:quad}))}\\
&=&2v.
\end{eqnarray*}
Similarly, $t_j\cdot t_k+t_k\cdot t_j=2\delta_{jk}v$ for all $j,k\in J$;  and
$s_j\cdot t_k+t_k\cdot s_j=0$ for all $j\in J-\{1\},k\in J$.

Next we consider the case that $u_0$ is present. If $j\ne 1$,
\begin{eqnarray*}
s_j\cdot u_0+u_0\cdot s_j&=&s_jv^*u_0+u_0v^*s_j=2\tp{s_j}{v}{u_0}\\
&=&-2i\tpc{u_j+\t{u}_j}{u_1+\t{u}_1}{u_0}\\
&=&-2i[\tp{u_j}{u_1}{u_0}+\tp{u_j}{\t{u}_1}{u_0}
+\tp{\t{u}_j}{u_1}{u_0}+\tp{\t{u}_j}{\t{u}_1}{u_0}].
\end{eqnarray*}
By Peirce calculus $\tp{u_j}{u_1}{u_0}$ is orthogonal to $u_1$ and
is a multiple of $u_j$, hence it is zero. Similarly, each of the other
three terms is zero, and 
similarly $t_j\cdot u_0+u_0\cdot t_j=0$ for all $j\in J$.
Finally $u_0\cdot u_0=-i\tpc{u_0}{u_1+\t{u}_1}{u_0}=v$ by (\ref{eq:govern}).
\qed 

\medskip 

As $Y$ is completely isometric to $Y_2(v)^{(v)}$, this
 completes the proof of Proposition~\ref{prop:1} and Theorem~\ref{thm:3}
in the case that
$Y$ is triple isomorphic to a Cartan factor of type 4.

\medskip

\subsection{Cartan factors of type 3}

We next prove Proposition~\ref{prop:1} and Theorem~\ref{thm:3}
in the case that 
$Y$ is triple isomorphic to the Cartan factor $S(H,J)$, of type 3. Again,
we let $A$ denote any von Neumann algebra containing $Y$.  

Let us recall that an {\it hermitian grid} (cf. \cite[p.\ 308]{DanFri87})
is a family $\{u_{ij}:i,j\in I\}$ of 
partial isometries satisfying
$u_{ij}=u_{ji}$; $u_{ij}\perp u_{kl}$ if $\{i,j\}\cap\{k,l\}=\emptyset$; 
$u_{ij}\vdash u_{ii}$ if $i\ne j$;
$u_{ij}\top u_{jk}$ if $i,j,k$ distinct; 
$\tp{u_{ij}}{u_{jk}}{u_{kl}}=u_{il}/2$
and $\tp{u_{ij}}{u_{jk}}{u_{ki}}=u_{ii}$
if $i\ne l$
and at least two of these elements are distinct;
 and all other
triple products are 0.

Let $\{u_{ij}:i,j\in \Lambda \}$
be an hermitian grid which is w*-total in $Y$ and let
$v$ denote the partial isometry  $\sum_iu_{ii}$ (the sum is w*-convergent since $u_{ii} \perp u_{jj}$) and note that it has
the property that $Y=Y_2(v)$. Let $\psi:
Y\rightarrow S(H,J)$ be the triple isomorphism determined by
$\psi(u_{ij})=U_{ij}$, where $\{U_{ij}\}$ denotes the canonical 
hermitian grid for $S(H,J)$, that is, $U_{ij}=\phi_{j} \otimes \phi_{i}+\phi_{i} \otimes \phi_{j}$ for $i\ne j$ and $U_{ii} = \phi_{i}
\otimes \phi_{i}$ for an orthonormal basis $\{ \phi_{\lambda} \}$ of $H$.

Note that isomorphisms of $JW^*$-triples (being isometries on spaces with
unique preduals) are automatically
w*-continuous.  Hence $\psi(v)=Id_H$ and
$\psi(u_{ij}^\sharp)=\tpc{\psi(v)}{\psi(u_{ij})}{\psi(v)}=
\psi(u_{ij})^*=U_{ij}^*=U_{ij}=\psi(u_{ij})$, so that $u_{ij}$ is self-adjoint
in $A_2(v)$. Here, $a^\sharp=va^*v$ denotes the involution in
 $A_2(v)^{(v)}$. 
Also, recall
 that the ternary product is the same whether it is computed in $A$
or in $A_2(v)^{(v)}$, that is, $xy^*z=xv^*(vy^*v)v^*z$.

Now define $e_{ij}=u_{ii}\cdot u_{ij}$, 
where we use $a\cdot b$ to denote 
the associative product in $A_2(v)^{(v)}$, that is,
$a\cdot b=av^*b$.

\begin{lem}
The collection $\{e_{ij}\}$ 
forms a system of matrix units in $A_2(v)^{(v)}$, that is,
\begin{description}
\item[(a)] $e_{ij}^\sharp=e_{ji}$, $e_{ij}\cdot e_{kl}=\delta_{jk}e_{il}$,
$v=\sum e_{ii}$.

Moreover,
\item[(b)] $u_{ii}\cdot u_{ij}=u_{ij}\cdot u_{jj}$ and $u_{ij}=e_{ij}+e_{ji}$.
\item[(c)] $\psi$ extends to a $^*$-isomorphism $\t{\psi}:
\mbox{sp}\, _{\CC}
\{e_{ij}\} \rightarrow \mbox{sp}\, _{\CC}\{E_{ij}\}$ satisfying $\t{\psi}(e_{ij})=E_{ij}$, where $E_{ij}=\phi_{j} \otimes
\phi_{i}$.
\end{description}
\end{lem}
\pf\
We first show these three identities:
\begin{equation}\label{eq:a1}
(u_{ii}\cdot u_{ij}-u_{ij}\cdot u_{jj})^\sharp 
(u_{ii}\cdot u_{ij}-u_{ij}\cdot u_{jj})=0,
\end{equation}
\begin{equation}\label{eq:a3}
(e_{ij}\cdot e_{kl})\cdot (e_{ij}\cdot e_{kl})^\sharp=0 \mbox{ for }j\ne k,
\end{equation}
\begin{equation}\label{eq:a2}
(e_{ij}\cdot e_{jl}-e_{il})\cdot (e_{ij}\cdot e_{jl}-e_{il})^\sharp=0.
\end{equation}

Note that $u_{ii}\cdot u_{ii}=u_{ii}v^*u_{ii}=u_{ii}(\sum_j u_{jj}^*)
u_{ii}=u_{ii}$, and that similarly, for $i\ne j$, $u_{ij}\cdot u_{ij}=
u_{ii}+u_{jj}$, $u_{ii}\cdot u_{jj}=0$, and
 $u_{ij}\cdot u_{ii}\cdot u_{ij}=u_{jj}$. 
Therefore, for $i\ne j$,
\begin{eqnarray*}
\lefteqn{
(u_{ii}\cdot u_{ij}-u_{ij}\cdot u_{jj})^\sharp\cdot(u_{ii}\cdot
u_{ij}-u_{ij}\cdot u_{jj})=}&\\
&=&(u_{ij}\cdot u_{ii}-u_{jj}\cdot u_{ij})\cdot (u_{ii}\cdot u_{ij}-u_{ij}
\cdot u_{jj})\\
&=&u_{ij}\cdot (u_{ii}\cdot u_{ii})\cdot u_{ij}-u_{jj}\cdot 
(u_{ij}\cdot u_{ii}\cdot u_{ij})\\
&-&(u_{ij}\cdot u_{ii}\cdot u_{ij})\cdot 
u_{jj}+u_{jj}\cdot (u_{ij}\cdot u_{ij})\cdot u_{jj}\\
&=&u_{ij}\cdot u_{ii}\cdot u_{ij}-u_{jj}\cdot 
u_{jj}-u_{jj}\cdot u_{jj}+u_{jj}\cdot (u_{ii}+u_{jj})\cdot u_{jj}\\
&=&u_{jj}-u_{jj}-u_{jj}+u_{jj}=0, \mbox{ proving (\ref{eq:a1}) and the 
first statement in (b)}.
\end{eqnarray*}

Next,  if $j\ne k$, then
\begin{eqnarray*}
(e_{ij}\cdot e_{kl})\cdot (e_{ij}\cdot e_{kl})^\sharp&=&
e_{ij}\cdot 
e_{kl}\cdot e_{lk}\cdot e_{ji}\\ 
&=&u_{ij}\cdot u_{jj}\cdot (u_{kl}\cdot u_{ll}\cdot u_{lk})\cdot 
u_{kk}\cdot u_{ji}\cdot u_{ii}\mbox{ (by (b))}\\
&=&u_{ij}\cdot u_{jj}\cdot (u_{kk}\cdot u_{kk})u_{ij}\cdot u_{ii}\\
&=&u_{ij}\cdot u_{jj}\cdot u_{kk}\cdot u_{ij}\cdot u_{ii}=0,
\mbox{ proving (\ref{eq:a3})}.
\end{eqnarray*}

Next,
\begin{eqnarray*}\lefteqn{ 
(e_{ij}\cdot e_{jl}-e_{il})\cdot (e_{ij}\cdot e_{jl}-e_{il})^\sharp=}&\\
&=&(e_{ij}\cdot e_{jl}-e_{il})\cdot (e_{lj}\cdot e_{ji}-e_{li})\\
&=&(u_{ij}\cdot u_{jj}\cdot u_{jl}\cdot u_{ll}-u_{il}
\cdot u_{ll})\cdot (u_{ll}\cdot u_{jl}\cdot u_{jj}\cdot 
u_{ij}-u_{ll}\cdot u_{li})\\
&=&-u_{ij}\cdot u_{jj}\cdot u_{jl}\cdot u_{ll}\cdot u_{li}+u_{il}
\cdot u_{ll}\cdot u_{li}\\
&+&u_{ij}\cdot u_{jj}\cdot (u_{jl}\cdot 
u_{ll}\cdot u_{jl})\cdot u_{jj}\cdot u_{ij}-
u_{il}\cdot u_{ll}\cdot u_{jl}\cdot u_{jj}\cdot u_{ij}\\
&=&-A+u_{ii}+u_{ij}\cdot u_{jj}\cdot u_{jj}\cdot u_{jj}\cdot u_{ij}-B\\
&=&-A+u_{ii}+u_{ii}-B,
\end{eqnarray*}
where $A=u_{ij}\cdot u_{jj}\cdot u_{jl}\cdot 
u_{ll}\cdot u_{li}$ and $B=u_{il}\cdot u_{ll}\cdot u_{jl}\cdot u_{jj}
\cdot u_{ij}$.

To prove (\ref{eq:a2}), it remains to show that $A=B=u_{ii}$. 
Here we need to distinguish cases.  Suppose first that
$i,j$ and $l$ are distinct.
  Then 
$\tp{u_{jl}}{u_{ll}}{u_{il}}=u_{ij}/2$ so that 
\begin{eqnarray*}
A&=&u_{ij}\cdot u_{jj}\cdot (2\tp{u_{jl}}{u_{ll}}{u_{il}}-u_{il}\cdot 
u_{ll}\cdot u_{jl})\\
&=&u_{ij}\cdot u_{jj}\cdot u_{ij}-u_{ij}\cdot u_{jj}\cdot (u_{il}
\cdot u_{ll})\cdot u_{jl}\\
&=&u_{ii}-u_{ij}\cdot u_{jj}\cdot (u_{ii}\cdot u_{il})\cdot u_{jl}
\mbox{ (by the first statement in (b)) }\\
&=&u_{ii}\mbox{ as required}.
\end{eqnarray*}

Also $B=u_{il}\cdot u_{ll}\cdot (2\tp{u_{jl}}{u_{jj}}{u_{ij}}-u_{ij}
\cdot u_{jj}\cdot u_{jl})=
u_{il}\cdot u_{ll}\cdot u_{li}-u_{il}\cdot u_{ll}\cdot u_{ij}
\cdot u_{jj}\cdot u_{jl}=u_{ii}$.

Now if $i=j$, then $A=u_{ii}\cdot u_{ii}
\cdot u_{il}\cdot u_{ll}\cdot u_{il}=u_{ii}^3=u_{ii}$ and $B=
u_{il}\cdot u_{ll}\cdot u_{il}\cdot 
u_{ii}\cdot u_{ii}=u_{ii}^3=u_{ii}$.  Similarly if $l=j$
or $i=l$, proving (\ref{eq:a2}).

Finally, since $u_{ij}u_{kk}=0$ if $k\not\in\{i,j\}$,
$e_{ij}+e_{ji}=u_{ij}\cdot u_{jj}+u_{ij}
\cdot u_{ii}=u_{ij}\cdot (\sum u_{kk})=u_{ij}$.  This completes the proof
of (a) and (b).

By the first statement in (b), we have
\[
e_{ij}^\sharp=(u_{ii}\cdot u_{ij})^\sharp=u_{ij}\cdot
u_{ii}=u_{jj}u_{ji}=e_{ji}.
\]

Since the system of matrix units $\{e_{ij}\}$ are linearly independent,
 $\tilde{\psi}$ defines a linear
isomorphism of $\mbox{sp}\, _{\CC}\{e_{ij}\}$ onto
$\mbox{sp}\, _{\CC}\{E_{ij}\}$ 
which is by construction a $^*$-isomorphism, proving
(c).\qed

\medskip

Clearly $\tilde{\psi}$ extends to a $^*$-isomorphism from the 
$C^*$-subalgebra $\overline{\mbox{sp}\, _{\CC}\{e_{ij}\}}^{\| \cdot
\|}$ of $A_2(v)^{(v)}$ onto the $C^*$-algebra of 
compact operators $K(H)$. By \cite[Lemma 1.14]{DanFri87}, 
$\tilde{\psi}$ extends to a w*-continuous 
isometry and,
hence, $^*$-isomorphism from
$\overline{\mbox{sp}\, _{\CC}\{e_{ij}\}}^{\mbox{w*}}$ onto $B(H)$. Since a
$^*$-isomorphism is completely isometric,  
the proof of Proposition~\ref{prop:1} 
and Theorem~\ref{thm:3} is completed
 in the case that
$Y$ is triple isomorphic to a Cartan factor of type 3.

\medskip

\section{Cartan factors of type 2}

In this section, we
 prove Proposition~\ref{prop:1} and Theorem~\ref{thm:3}
in the case that 
$Y$ is triple isomorphic to the Cartan factor $A(H,J)$ of type 2. Again,
$A$ denotes any von Neumann algebra containing $Y$. Since 
$A(\CC^3,J)$ is triple isomorphic to $B(\CC,\CC^3)$,
which is covered in section 7, and $A(\CC^4,J)$
is triple isomorphic to $\mbox{Sp}\, (6)$, which was covered in section
3, we may and shall assume that $\dim H > 4$. 

Let us recall (\cite[p.\ 317]{DanFri87}
that a {\it symplectic grid} is a family $\{u_{ij}:i,j\in I,i\ne j\}$ of 
minimal partial isometries satisfying
 $u_{ij}=-u_{ji}$; $u_{ij}\top u_{kl}$ if $\{i,j\}\cap \{k,l\}\ne\emptyset$;
$u_{ij}\perp u_{kl}$ if $\{i,j\}\cap \{k,l\}=\emptyset$; $2\tp{u_{ij}}{u_{il}}{u_{kl}}=u_{kj}$ for distinct
$i,j,k,l$; and all other triple products vanish. 
The fact that each $u_{ij}$ is  minimal can be
expressed by
\begin{equation}\label{eq:12bis}
u_{ij}u_{kl}^*u_{ij}=\delta_{(i,j),(k,l)}u_{ij}.
\end{equation}

Let $\{u_{ij}\}$
be a symplectic grid which is w*-total in $Y$.
Let $\psi:
Y\rightarrow A(H,J)$ be the triple isomorphism determined by
$\psi(u_{ij})=U_{ij}$, where $\{U_{ij}\} = \phi_{j} \otimes \phi_{i} - \phi{i} \otimes \phi_{j}$ for an orthonormal
basis $\{ \phi_{\lambda} \}$ of $H$.

\begin{lem}\label{lem:4.1}
For any indices $i,j,k,l,m$,
\begin{equation}\label{eq:620}
u_{ik}u_{kl}^*u_{il}=u_{ij}u_{jm}^*u_{im},
\end{equation}
 and for $1\le i\le n$,
the elements
$e_{ii}$  unambiguously defined by $e_{ii}=u_{ij}u_{jm}^*u_{im}$ are non-zero 
orthogonal partial isometries in $A$.
\end{lem}
 \pf\
We shall use repeatedly the fact that $u_{ij}=-u_{ji}$. 

Suppose that $i,j,k,l$ are distinct. Then $u_{ij}u_{kl}^*=0$ and therefore
\begin{eqnarray}\label{eq:4251}\nonumber
u_{ik}u_{kl}^*u_{il}&=&2\tp{u_{ij}}{u_{ij}}{u_{ik}}u_{kl}^*u_{il}\\\nonumber
&=&(u_{ij}u_{ij}^*u_{ik}+u_{ik}u_{ij}^*u_{ij})u_{kl}^*u_{il}\\
&=&u_{ij}u_{ij}^*u_{ik}u_{kl}^*u_{il}+0\\\nonumber
&=&u_{ij}(u_{kl}u_{ik}^*u_{ij}+u_{ij}u_{ik}^*u_{kl})^*u_{il}\\\nonumber
&=&2u_{ij}\tp{u_{ij}}{u_{ik}}{u_{kl}}^*u_{il}\\\nonumber
&=&u_{ij}(-u_{lj})^*u_{il}=u_{ij}u_{jl}^*u_{il}.
\end{eqnarray}

Similarly, if $m,l,i,k$ are distinct, 
by replacing $u_{il}$ by $2\tp{u_{im}}{u_{im}}{u_{il}}$, we
obtain $u_{ik}u_{kl}^*u_{il}=u_{ik}u_{km}^*u_{im}$. Indeed,
\begin{eqnarray}\label{eq:4252}\nonumber
u_{ik}u_{kl}^*u_{il}&=&u_{ik}u_{kl}^*2\tp{u_{im}}{u_{im}}{u_{il}}\\\nonumber
&=&u_{ik}u_{kl}^*(u_{im}u_{im}^*u_{il}+u_{il}u_{im}^*u_{im})\\
&=&u_{ik}u_{kl}^*u_{il}u_{im}^*u_{im}\\\nonumber
&=&2u_{ik}\tp{u_{kl}}{u_{il}}{u_{im}}^*u_{im}\\\nonumber
&=&u_{ik}u_{km}^*u_{im}.
\end{eqnarray}

Together, (\ref{eq:4251}) and (\ref{eq:4252}) show (\ref{eq:620})
and thus $e_{ii}$ is well defined.

We next show that $e_{ii}\ne 0$.
Suppose instead that $e_{ii}=0$ for some $i$. For $i,k,l$ distinct,
$u_{ik}\top u_{kl}$ and $u_{kl}\top u_{il}$, so
\begin{eqnarray}\label{eq:*}\nonumber
u_{kl}&=&u_{ik}u_{ik}^*u_{kl}+u_{kl}u_{ik}^*u_{ik}\\\nonumber
&=&u_{ik}u_{ik}^*(u_{il}
u_{il}^*u_{kl}+u_{kl}u_{il}^*u_{il})+
(u_{il}u_{il}^*u_{kl}+u_{kl}u_{il}^*u_{il})u_{ik}^*u_{ik}\\
&=&(u_{ik}u_{ik}^*u_{il}u_{il}^*u_{kl}+u_{ik}e_{ii}^*u_{il})+(u_{il}e_{ii}^*u_{ik}+u_{kl}u_{il}^*u_{il})
u_{ik}^*u_{ik})\\\nonumber
&=&(u_{ik}u_{ik}^*u_{il}u_{il}^*u_{kl}+0)+(0+u_{kl}u_{il}^*u_{il})
u_{ik}^*u_{ik})\\\nonumber
&=& L_{ik}L_{il}u_{kl}+u_{kl}R_{il}R_{ik},
\end{eqnarray}
where $L_{ik}=u_{ik}u_{ik}^*$ and $R_{ik}=u_{ik}^*u_{ik}$ denote the
left and right support projections of $u_{ik}$.

By the definition of symplectic grid, if $p,k,l,m$ are distinct
(recall that $n\ge 5$), then
$u_{pm}=2\tp{u_{pk}}{u_{kl}}{u_{ml}}$.
However,  by (\ref{eq:*}) and the commutativity of
  the support projections associated with $u_{il}$ and $u_{ik}$ (see 
Lemma~\ref{lem:4.6}),
\begin{eqnarray*}
u_{pm}&=& 2\tp{u_{pk}}{u_{kl}}{u_{ml}}\\
&=&u_{pk}u_{kl}^*u_{ml}+u_{ml}u_{kl}^*u_{pk}\\
&=&
u_{pk}(L_{ik}L_{il}u_{kl}+u_{kl}R_{il}R_{ik})^*u_{ml}+
u_{ml}(L_{ik}L_{il}u_{kl}+u_{kl}R_{il}R_{ik})^*u_{nk}\\
&=&u_{pk}u_{kl}^*L_{il}(L_{ik}u_{ml})+(u_{pk}R_{il})R_{ik}u_{kl}^*u_{ml}\\
&&+u_{ml}u_{kl}^*L_{ik}(L_{il}u_{pk})
+(u_{ml}R_{ik})R_{il}u_{kl}^*u_{pk}=0,
\end{eqnarray*}
which is a contradiction.
Thus $e_{ii}\ne 0$ and it is a partial isometry:
\begin{eqnarray*}
e_{ii}e_{ii}^*e_{ii}&=&u_{ik}u_{kl}^*u_{il}u_{il}^*u_{kl}u_{ik}^*(u_{ik}
u_{kl}^*u_{il})\\
&=&-u_{ik}u_{kl}^*u_{il}u_{il}^*(u_{kl}u_{ik}^*u_{il})u_{kl}^*u_{ik}\\
&=&u_{ik}u_{kl}^*(u_{il}u_{il}^*u_{il})u_{ik}^*u_{kl}u_{kl}^*u_{ik}\\
&=&u_{ik}(u_{kl}^*u_{il}u_{ik}^*)u_{kl}u_{kl}^*u_{ik}\\
&=&-u_{ik}u_{ik}^*u_{il}(u_{kl}^*u_{kl}u_{kl}^*)u_{ik}\\
&=&-u_{ik}u_{ik}^*(u_{il}u_{kl}^*u_{ik})\\
&=&u_{ik}u_{ik}^*e_{ii}\\
&=&u_{ik}u_{ik}^*(u_{ik}u_{km}^*u_{im})\\
&=&u_{ik}u_{km}^*u_{im}=e_{ii}.
\end{eqnarray*}

Finally, to show orthogonality, take $i,j,l,p$ distinct
and note that 
\[
e_{ii}^*e_{jj}=(u_{il}u_{lk}^*u_{ki})^*u_{jp}u_{pm}^*u_{mj}=
u_{ki}^*u_{lk}u_{il}^*u_{jp}u_{pm}^*u_{mj}=0
\]
and similarly $e_{ii}e_{jj}^*=0$. \qed

\begin{lem}\label{lem:4.2}
With the above notation,
\begin{description}
\item[(a)] $e_{ii}u_{ij}^*e_{ii}=e_{ii}e_{jj}^*e_{ii}=0$ for $i\ne j$.
\item[(b)] $\{e_{ii}\}\cup \{u_{ij}\}$ is a linearly independent set.
\item[(c)] $u_{ij}\perp e_{kk}$ for $k\not\in\{i,j\}$, that is, $u_{ij}\in A_0(e_{kk})$.
\item[(d)] $\tp{e_{ii}}{e_{ii}}{u_{ij}}=u_{ij}/2=\tp{e_{jj}}{e_{jj}}{u_{ij}}$, that is,
$u_{ij}\in A_1(e_{ii})\cap A_1(e_{jj})$.
\end{description}
\end{lem}

\pf\ 
\begin{eqnarray*}
e_{ii}u_{ij}^*e_{ii}&=&u_{ik}u_{kl}^*u_{il}u_{ij}^*u_{ik}u_{kl}^*u_{il}=
-u_{il}u_{kl}^*u_{ik}u_{ij}^*u_{ik}u_{kl}^*u_{il}\\
&=&-u_{il}u_{kl}^*\tp{u_{ik}}{u_{ij}}{u_{ik}}u_{kl}^*u_{il}=0,
\end{eqnarray*}
since $\tp{u_{ik}}{u_{ij}}{u_{ik}}\in A_{2-1+2}(u_{ik})=\{0\}$
 by Peirce calculus. Also
\[
e_{ii}e_{jj}^*e_{ii}=u_{ik}u_{kl}^*(u_{il}u_{jp}^*)u_{pm}u_{jm}^*u_{ik}
u_{kl}^*u_{il}=0,
\]
 proving (a), and (b) follows
immediately from (a).

To prove (c), note first that $e_{kk}^*u_{ij}=(u_{kl}u_{lm}^*u_{mk})^*u_{ij}=
u_{mk}^*u_{lm}u_{kl}^*u_{ij}$ which is zero if $\{i,j\}\cap \{l,k\}
=\emptyset$, and similarly
$u_{ij}e_{kk}^*=0$. Thus $2\tp{u_{ij}}{u_{ij}}{e_{kk}}=
u_{ij}u_{ij}^*e_{kk}+e_{kk}u_{ij}^*u_{ij}=0$,
 which is equivalent to $u_{ij}\in A_0(e_{kk})$.

Finally, we shall
show  assertion (d): 
\begin{eqnarray*}
2\tp{e_{ii}}{e_{ii}}{u_{ij}}&=&e_{ii}e_{ii}^*u_{ij}+u_{ij}e_{ii}^*e_{ii}\\
&=&u_{ip}u_{pm}^*u_{im}u_{il}^*u_{kl}u_{ik}^*u_{ij}+
u_{ij}u_{ik}^*u_{lk}u_{il}^*u_{im}u_{mp}^*u_{ip}\\
&=&2u_{ip}\tp{u_{pm}}{u_{im}}{u_{il}}^*u_{kl}u_{ik}^*u_{ij}+
2u_{ij}u_{ik}^*u_{lk}\tp{u_{il}}{u_{im}}{u_{mp}}^*u_{ip}\\
&&(\mbox{since } u_{pm}^*u_{kl}=u_{kl}u_{pm}^*=0)\\
&=&u_{ip}u_{pl}^*u_{kl}u_{ik}^*u_{ij}-u_{ij}u_{ik}^*u_{lk}u_{pl}^*u_{ip}\\
&=&2u_{ip}\tp{u_{pl}}{u_{kl}}{u_{ik}}^*u_{ij}-2u_{ij}
\tp{u_{ik}}{u_{lk}}{u_{pl}}^*u_{ip}\\
&=&u_{ip}(-u_{pi})^*u_{ij}+u_{ij}u_{ip}^*u_{ip}\\
&=&2\tp{u_{ip}}{u_{ip}}{u_{ij}}=u_{ij}.\qed
\end{eqnarray*}

%

\begin{lem}
Define, for $i\ne j$, $e_{ij}=e_{ii}e_{ii}^*u_{ij}e_{jj}^*e_{jj}$ 
{\rm (}product in $A${\rm )}. 
Then $u_{ij}=e_{ij}-e_{ji}$ and
\begin{description}
\item[(a)] $\{e_{ij}\}$ is a system of matrix units in the \csa\ 
$A_2(v)^{(v)}$. 
\item[(b)] $\psi$ extends to a $^*$-isomorphism $\t{\psi}:
\mbox{sp}\, _{\CC}\{e_{ij}\} \rightarrow \mbox{sp}\, _{\CC}\{E_{ij}\}$ satisfying $\t{\psi}(e_{ij})=E_{ij}$, where
$\{E_{ij}\} = \phi_{j} \otimes \phi_{i}$.
\end{description}
\end{lem}
\pf\
By definition, $e_{ji}=e_{jj}e_{jj}^*u_{ji}e_{ii}^*e_{ii}=-e_{jj}e_{jj}^*
u_{ij}e_{ii}^*e_{ii}=-e_{ij}$, and by
Lemmas~\ref{lem:4.1} and ~\ref{lem:4.2}, 
\begin{eqnarray*}
e_{ij}-e_{ji}&=&   
e_{ii}e_{ii}^*u_{ij}e_{jj}^*e_{jj}+e_{jj}e_{jj}^*u_{ij}
e_{ii}^*e_{ii}\\
&=&2\tp{e_{ii}}{e_{ii}}{u_{ij}}e_{jj}^*e_{jj}+e_{jj}e_{jj}^*
2\tp{e_{ii}}{e_{ii}}{u_{ij}}\\
&=&u_{ij}e_{jj}^*e_{jj}+e_{jj}e_{jj}^*u_{ij}\\
&=&2\tp{e_{jj}}{e_{jj}}{u_{ij}}=u_{ij}.
\end{eqnarray*}
Since $v=\sum e_{kk}$, to prove (a)
it remains to  show that $e_{ij}v^*e_{lk}=\delta_{jl}e_{ik}$ 
and $ve_{ij}^*v=e_{ji}$.

In the first place, if 
$j\ne l$, then $e_{ij}v^*e_{lk}=e_{ii}e_{ii}^*u_{ij}
e_{jj}^*(e_{jj}v^*e_{ll})e_{ll}^*u_{lk}e_{kk}^*e_{kk}=0$.

Now consider  the case $j=l$, so that 
\[
e_{ij}v^*e_{jk}=e_{ii}e_{ii}^*u_{ij}e_{jj}^*e_{jj}(\sum_q e_{qq})^*
e_{jj}e_{jj}^*u_{jk}e_{kk}^*e_{kk}
=e_{ii}e_{ii}^*u_{ij}e_{jj}^*u_{jk}e_{kk}^*e_{kk}.
\]
 There are five
cases to prove:
\begin{itemize}
\item $e_{ii}v^*e_{ii}=e_{ii}$; this is true since
$e_{ii}v^*e_{ii}=e_{ii}e_{ii}^*e_{ii}=e_{ii}$.
\item $e_{ii}v^*e_{ik}=e_{ik}$, $k\ne i$; since $e_{ii}v^*e_{ik}=e_{ii}v^*
e_{ii}e_{ii}^*u_{ik}e_{kk}^*e_{kk}=e_{ii}e_{ii}^*u_{ik}e_{kk}^*e_{kk}=e_{ik}$.
\item $e_{ij}v^*e_{jj}=e_{ij}$, $i\ne j$; since 
$e_{ij}v^*e_{jj}=e_{ii}e_{ii}^*
u_{ij}e_{jj}^*e_{jj}v^*e_{jj}=e_{ij}$.
\item $e_{ij}v^*e_{ji}=e_{ii}$, $i\ne j$; since
\begin{eqnarray*}
e_{ij}v^*e_{ji}&=&e_{ii}e_{ii}^*u_{ij}e_{jj}^*e_{jj}v^*e_{jj}e_{jj}^*u_{ji}
e_{ii}^*e_{ii}\\
&=&e_{ii}e_{ii}^*u_{ij}e_{jj}^*u_{ji}
e_{ii}^*e_{ii}\\
&=&e_{ii}e_{ii}^*u_{ij}u_{jk}^*u_{kl}u_{jl}^*u_{ji}
e_{ii}^*e_{ii}\\
&=&2e_{ii}e_{ii}^*\tp{u_{ij}}{u_{jk}}{u_{kl}}u_{jl}^*u_{ji}
e_{ii}^*e_{ii}\\
&&(\mbox{by Lemma~\ref{lem:4.2}(c), }e_{ii}^*u_{kl}=0)\\
&=&-e_{ii}e_{ii}^*u_{il}u_{jl}^*u_{ji}
e_{ii}^*e_{ii}\\
&=&e_{ii}e_{ii}^*e_{ii}
e_{ii}^*e_{ii}=e_{ii}.
\end{eqnarray*}
\item $e_{ij}v^*e_{jk}=e_{ik}$, $i,j,k$ distinct; since
\begin{eqnarray*}
e_{ij}v^*e_{jk}&=&e_{ii}e_{ii}^*u_{ij}(e_{jj}^*e_{jj}e_{jj}^*e_{jj}
e_{jj}^*)u_{jk}e_{kk}^*e_{kk}\\
&=&e_{ii}e_{ii}^*u_{ij}e_{jj}^*u_{jk}e_{kk}^*e_{kk}\\
&=&e_{ii}e_{ii}^*u_{ij}(u_{jm}u_{mp}^*u_{jp})^*u_{jk}e_{kk}^*e_{kk}\\
&=&e_{ii}e_{ii}^*u_{ij}u_{jp}^*u_{mp}u_{jm}^*u_{jk}e_{kk}^*e_{kk}\\
&=&2e_{ii}e_{ii}^*\tp{u_{ij}}{u_{jp}}{u_{mp}}
u_{jm}^*u_{jk}e_{kk}^*e_{kk}\\
&=&e_{ii}e_{ii}^*u_{im}u_{jm}^*u_{jk}e_{kk}^*e_{kk}\\
&=&2e_{ii}e_{ii}^*\tp{u_{im}}{u_{jm}}{u_{jk}}
e_{kk}^*e_{kk}\\
&=&e_{ii}e_{ii}^*u_{ik}e_{kk}^*e_{kk}=e_{ik}.
\end{eqnarray*}
\end{itemize}

Finally, 
\begin{eqnarray*}
e_{ij}&=&e_{ii}e_{ii}^*u_{ij}e_{jj}^*e_{jj}\\
&=&u_{ik}u_{kl}^*u_{il}u_{il}^*u_{kl}(u_{ik}^*u_{ij}
u_{jm}^*)u_{pm}u_{jp}^*u_{jp}u_{pm}^*u_{jm}\\
&=&2u_{ik}u_{kl}^*u_{il}u_{il}^*u_{kl}
\tp{u_{ik}}{u_{ij}}{u_{jm}}^*
u_{pm}u_{jp}^*u_{jp}u_{pm}^*u_{jm}\\
&=&-u_{ik}u_{kl}^*u_{il}u_{il}^*u_{kl}(u_{mk}^*
u_{pm}u_{jp}^*)u_{jp}u_{pm}^*u_{jm}\\
&=&u_{ik}u_{kl}^*u_{il}(u_{il}^*u_{kl}u_{jk}^*)u_{jp}u_{pm}^*u_{jm}\\
&=&-u_{ik}u_{kl}^*u_{il}u_{ij}^*u_{jp}u_{pm}^*u_{jm}\\
&=&-e_{ii}u_{ij}^*e_{jj}
\end{eqnarray*}
and therefore $ve_{ij}^*v=-ve_{jj}^*u_{ij}e_{ii}^*v=-e_{jj}e_{jj}^*u_{ij}
e_{ii}^*e_{ii}=
e_{jj}e_{jj}^*u_{ji}
e_{ii}^*e_{ii}=e_{ji}$.
This completes the proof of (a).

Since the system of matrix units $\{e_{ij}\}$ are linearly independent,
 $\tilde{\psi}$ defines a $^*$-isomorphism 
of $\mbox{sp}\, _{\CC}\{e_{ij}\}$ onto 
$\mbox{sp}\, _{\CC}\{E_{ij}\}$, proving
(b).\qed

\medskip

By the same method used in Section 3 for the Type 3 case, $\tilde{\psi}$ extends to a $^*$-isomorphism of
$\overline{\mbox{sp}\, _{\CC}\{e_{ij}\}}^{\mbox{w*}}$ 
onto $B(H)$. The proof of
Proposition~\ref{prop:1} and Theorem~\ref{thm:3}
is thus complete in the case that
$Y$ is triple isomorphic to a Cartan factor of type 2.

\section{Cartan factors of type 1}

In this and the next section, we
 prove Proposition~\ref{prop:1} and Theorem~\ref{thm:3}
in the case that 
$Y$ is triple isomorphic to a Cartan factor of type 1. 
This turns out to be more complicated than the other types, especially
in the case that $Y$ is of rank 1 (Hilbertian).  Except for some
important preliminary cases (see subsection \ref{5.3}), 
the rank 1 case is proved in section 7.

Let $\{u_{ij}:i\in \Lambda,j\in \Sigma\}$
be a {\it rectangular grid} which is w*-total in $Y$.
Recall (\cite[p.\ 313]{DanFri87}
that this means that each $u_{ij}$ is a minimal partial
isometry,
$u_{jk}\perp u_{il}$ if $i\ne j$ and $k\ne l$; $u_{jk}\top
u_{il}$ if either $j=i,k\ne l$ or $j\ne i,k=l$; 
\begin{equation}\label{eq:703}
\tp{u_{jk}}{u_{jl}}{u_{il}}=u_{ik}/2\mbox{ if }j\ne i\mbox{ and }k\ne l;
\end{equation}
 and all
other triple products are zero.

We shall assume throughout this section that $Y$ is triple isomorphic to
$B(H,K)$, that is, $|\Lambda|=\dim K$ and $|\Sigma|=\dim H$.  
Specifically, by \cite[p.317 and Lemma 1.14]{DanFri87}, this
means that the map $u_{ij}\mapsto E_{ij}$ extends to a triple isomorphism of $Y$ onto
$B(H,K)$, where $E_{ij}=\phi_{i} \otimes \psi_{j}$ for orthonormal bases 
$\{ \psi_{j}:j\in\Sigma \}$ in $H$ and $\{ \phi_{i}:i\in \Lambda \}$ in $K$.

\subsection{A special case}

Note that the canonical rectangular grid $\{ E_{ij} \}$ for $B(H,K)$ satisfies
 $E_{ij}E_{ik}^*=\ip{\psi_k}{\psi_j}\phi_i\otimes\phi_i=0$ 
for $j\ne k$ and all $i$;
and
$E_{ik}^*E_{jk}=0$ for $i\ne j$ and all $k$.

\begin{lem}\label{lem:4.5} With $Y$ as above, assume that
 for some fixed values of $i\in\Lambda,k,l\in\Sigma$,
$u_{il}u_{ik}^*=0$ and $k\ne l$, or for some fixed values of $i,j\in\Lambda,
k\in\Sigma$,
$u_{ik}^*u_{jk}=0$ and $i\ne j$. 
  Then 
\begin{description}
\item[(a)] for
all $j\in\Lambda,p,q\in\Sigma$ with $p\ne q$,
$u_{jp}u_{jq}^*=0$ and for all $p,q\in\Lambda,r\in\Sigma$ with $p\ne q$,
$u_{pr}^*u_{qr}=0$.
\item[(b)] $Y$ is a ternary subtriple of $A$ which is
 ternary isomorphic and completely 
isometric to $B(H,K)$.
\end{description}
\end{lem}
\pf\
We shall give the proof in the case that $u_{il}u_{ik}^*=0$. The other
case follows by symmetry. 

We first take care of the ``$i^{\mbox{th}}$-row,'' where
$i,k,l$ are the fixed values. If $p\not\in\{k,l\}$,
\begin{eqnarray}\label{eq:33}\nonumber
u_{ip}u_{ik}^*&=&2\tp{u_{il}}{u_{il}}{u_{ip}}u_{ik}^*\\\nonumber
&=&(u_{il}u_{il}^*u_{ip}+u_{ip}u_{il}^*u_{il})u_{ik}^*\\
&=&u_{il}(u_{il}^*u_{ip}u_{ik}^*)\mbox{ (by assumption)}\\\nonumber
&=&-u_{il}u_{ik}^*u_{ip}u_{il}^*\mbox{ (since }
\tp{u_{ik}}{u_{ip}}{u_{il}}=0)\\\nonumber
&=&0.
\end{eqnarray}

Thus, if $q\not\in\{p,k\}$,
\begin{eqnarray}\label{eq:44}\nonumber
u_{ip}u_{iq}^*&=&2u_{ip}\tp{u_{ik}}{u_{ik}}{u_{iq}}^*\\
&=&u_{ip}(u_{ik}u_{ik}^*u_{iq}+u_{iq}u_{ik}^*u_{ik})^*\\\nonumber
&=&(u_{ip}u_{iq}^*u_{ik})u_{ik}^*\mbox{ (by (\ref{eq:33}))}\\\nonumber
&=&-u_{ik}u_{iq}^*u_{ip}u_{ik}^*=0.
\end{eqnarray}

This proves the first statement in (a) when $j$ has the value $i$.
We next take care of the ``$j^{\mbox{th}}$-row'' (if it exists). If $p\ne q$
and $j\ne i$,
\begin{eqnarray}\label{eq:jthrow}\nonumber
u_{jp}u_{jq}^*&=&2\tp{u_{ip}}{u_{iq}}{u_{jq}}u_{jq}^*\\\nonumber
&=&(u_{ip}u_{iq}^*u_{jq}+u_{jq}u_{iq}^*u_{ip})u_{jq}^*\\
&=&u_{jq}u_{iq}^*u_{ip}u_{jq}^*\mbox{ (by (\ref{eq:44}))}\\\nonumber
&=&2u_{jq}u_{iq}^*u_{ip}\tp{u_{iq}}{u_{ip}}{u_{jp}}^*\\\nonumber
&=&u_{jq}u_{iq}^*u_{ip}(u_{iq}^*u_{ip}u_{jp}^*+u_{jp}^*u_{ip}u_{iq}^*)\\
\nonumber
&=&0 \mbox{ (by (\ref{eq:44}) and the minimality of $u_{ip}$)}.
\end{eqnarray}

This proves the first statement in (a).
We complete the proof of (a), by taking
 care of the ``columns,'' (if they exist).  
For all $p,q,r$ with $p\ne q$, choose $s\ne r$.  Then
\begin{eqnarray*}
u_{pr}^*u_{qr}&=&2u_{pr}^*\tp{u_{qs}}{u_{ps}}{u_{pr}}=u_{pr}^*
(u_{pr}u_{ps}^*u_{qs}+u_{qs}u_{ps}^*u_{pr})\\
&=&u_{pr}^*(u_{pr}u_{ps}^*)u_{qs}  
+(u_{pr}^*u_{qs})u_{ps}^*u_{pr}=0, 
\end{eqnarray*}
 by orthogonality of $u_{qs}$ and $u_{pr}$ for $s \ne r$ and
by (\ref{eq:jthrow}). 

By (a), and the separate w*-continuity of multiplication, $Y$ is a ternary subtriple (closed under $(a,b,c)\mapsto
ab^*c$). Furthermore,  $\mbox{sp}_{\CC}\{u_{ij}\}$ is ternary
isomorphic and isometric to $\mbox{sp}_{\CC}\{E_{ij}\}$ 
via $u_{ij}\mapsto E_{ij}$ and we can again use \cite[Lemma 1.14]{DanFri87} to extend the
map to a ternary isomorphism of $Y$ onto $B(H,K)$, which then is a complete isometry. 
 \qed

\medskip 


By the same arguments, we also have the following.
\begin{lem}\label{lem:4.5prime} With $Y$ as above, assume
that for some fixed values of $i,k,l$,
$u_{il}^*u_{ik}=0$ and $k\ne l$, or $u_{ik}u_{jk}^*=0$ and $i\ne j$.
  Then 
\begin{description}
\item[(a)] for
all $i,k,l$ with $k\ne l$,
$u_{il}^*u_{ik}=0$ and for all $i,j,k$ with $i\ne j$,
$u_{ik}u_{jk}^*=0$.
\item[(b)] $Y$ is ternary isomorphic and completely
isometric to $B(K,H)$.
\end{description}
\end{lem}

It is convenient to single out the rank one case.

\begin{cor}\label{cor:1}
Let $Y$ be triple isomorphic to a Cartan factor of type 1 and rank~ 1 and
denote by $\{u_{\lambda}\}$ a rank 1 rectangular grid for $Y$.
\begin{description}
\item[(a)] If 
$u_iu_j^*=0$ for some $i\ne j$, then $Y$ is completely isometric to
$B(H,\CC)$.
\item[(b)] If 
$u_i^*u_j=0$ for some $i\ne j$, then $Y$ is completely isometric to
$B(\CC,K)$.
\end{description}
\end{cor}

\subsection{The case of rank 2 or more}\label{sect:5.2}

The following simple lemma will be useful in this and the next
section. Part (b) of it is referred to
as ``hopping''.
\begin{lem}\label{lem:4.6}
Let $u,v,w$ be partial isometries.
\begin{description}
\item[(a)] If $u$ and $w$ are colinear, then the support projections $uu^*,ww^*$ commute as do
$u^*u,w^*w$.
 \item[(b)] If $v$ and $w$ are each colinear with $u$, then
$uu^*vw^*=vw^*uu^*$ and $u^*uv^*w=v^*wu^*u$.
\end{description}
\end{lem}
\pf\
We prove (b) first. Since $uu^*v+vu^*u=v$ and $uu^*w+wu^*u=w$,
$
(uu^*v)w^*=(v-vu^*u)w^*=vw^*-v(u^*uw^*)=vw^*-v(w^*-w^*uu^*)=
vw^*uu^*$.
Similarly for the second statement.  To prove (a) use the same argument:
$uu^*ww^*
=(w-wu^*u)w^*=ww^*-w(u^*uw^*)=ww^*-w(w^*-w^*uu^*)=ww^*uu^*$.\qed 

\medskip

Justified by Lemmas~\ref{lem:4.5} and \ref{lem:4.5prime}
we may now assume in the rest of this subsection~\ref{sect:5.2}, 
without loss of generality,
 that 
\begin{equation}\label{eq:466}
u_{ik}u_{ij}^*\ne 0\mbox{ and }u_{ik}^*u_{ij}\ne 0\mbox{ for all $i\in\Lambda,
j,k\in\Sigma$},
\end{equation}
and
\begin{equation}\label{eq:15bis}
u_{ik}u_{jk}^*\ne 0\mbox{ and }u_{ik}^*u_{jk}\ne 0\mbox{ for all $i,j\in
\Lambda,k\in\Sigma$}.
\end{equation}

\begin{lem}\label{lem:5.4}
Suppose that $Y$ is triple isomorphic to a Cartan factor $B(H,K)$
of type 1 and rank
at least 2, let $\{u_{ij}:i\in \Lambda,j\in \Sigma\}$ 
be a rectangular grid for $Y$ and suppose that
{\rm (\ref{eq:466})} and 
{\rm (\ref{eq:15bis})} hold. 

Then for all $i\in\Lambda$ and $j\in\Sigma$, the projections
\[
L_i:=\prod_{k \in \Sigma} u_{ik}u_{ik}^*
\mbox{ and }
R_j:=\prod_{l \in \Lambda} u_{lj}^*u_{lj} 
\]
are non-zero.
\end{lem}
\pf\
Note that by Lemma~\ref{lem:4.6}, the above are products of commuting
projections. We shall show that $L_i\ne 0$, the proof for $R_j$ being
similar.

Suppose the assertion is false, that is, for some $i\in\Lambda$,
\[
\prod_{k \in \Sigma}u_{ik}u_{ik}^*=0.
\]
Choose a finite subset $S \subseteq \Sigma$ and denote it by
$\{1,2, \cdots n\}$. Choose a $j \ne i$ and an 
$l\in S-\{1\}$.  Since $u_{il}\top u_{i1}$,

\begin{eqnarray*}
u_{il}u^{\ast}_{jl}&=&
u_{i1}u_{i1}^{\ast}u_{il}u^{\ast}_{jl} + 
u_{il}u^{\ast}_{i1}u_{i1}u^{\ast}_{jl}\\
&=& u_{i1}u_{i1}^{\ast}u_{il}u^{\ast}_{jl}\\
&=&u_{i1}u_{i1}^*u_{i2}u_{i2}^*u_{il}u_{jl}^*\\
&=& \cdots\\
&=&\left(\prod_{k=1,k\ne l}^n u_{ik}u_{ik}^*\right)u_{il}u^{\ast}_{jl}\\
&=&\left(\prod_{k=1}^n u_{ik}u_{ik}^*\right)u_{il}u^{\ast}_{jl} 
\end{eqnarray*}

Since $\prod_{k \in \Sigma}u_{ik}u_{ik}^*$ is the
 w*-limit of the net $\{ \prod_{k \in
S}u_{ik}u_{ik}^*\}_{|S| < \infty}$, it follows by separate w*-continuity of multiplication that $u_{il}u^{\ast}_{jl}=0$,
which contradicts (\ref{eq:15bis}). 
\qed

\begin{lem} \label{lem:5.5} Let $Y$ be as in Lemma~\ref{lem:5.4}.
Let \[
p=\sum_{i \in \Lambda}\prod_{k \in \Sigma} u_{ik}u_{ik}^*,
\]
which is a sum of non-zero orthogonal projections.
The maps $Y\ni y\mapsto py\in pY$ and $Y\ni y\mapsto (1-p)y\in (1-p)Y$
are completely contractive
triple isomorphisms. Also, $pY\perp(1-p)Y$.
\end{lem}
\pf\ 
We begin by showing
that $\{pu_{ij}\}$ is a rectangular grid which is w*-total in its
w*-closure.
We start by showing that $pu_{ij}$ is a minimal partial isometry
using the criterion (\ref{eq:12bis}).
We have
\[
pu_{ij}(pu_{kl})^*pu_{ij}=pu_{ij}u_{kl}^*pu_{ij}
=pu_{ij}u_{kl}^*\left(\sum_{q\in\Lambda} u_{q1}u_{q1}^*\cdots 
u_{qn}u_{qn}^*\right)u_{ij}.
\]
By Lemma~\ref{lem:4.6}(a), this is zero if $k\ne i$. For $k=i$, we have
\[
pu_{ij}(pu_{kl})^*pu_{ij}=pu_{ij}u_{il}^*u_{i1}u_{i1}^*\cdots 
u_{in}u_{in}^*u_{ij}.
\]
which is zero for $j\ne l$, since by Peirce calculus, 
$u_{ij}u_{il}^*u_{ij}=0$.  On the other hand
\[
(pu_{ij})(pu_{ij})^*(pu_{ij})=
pu_{ij}u_{ij}^*u_{i1}u_{i1}^*\cdots u_{in}u_{in}^*u_{ij}=pu_{ij},
\]
and it is non-zero by Lemma~\ref{lem:5.4}.
This proves that $pu_{ij}$ is a minimal partial isometry.

We next show that $pu_{jk}\perp pu_{il}$ for $i \ne j$ and $k \ne l$.
On the one hand, $pu_{jk}(pu_{il})^*=pu_{jk}u_{il}^*p=0$; and on the other 
hand, 
\[
(pu_{il})^*pu_{jk}=u_{il}^*pu_{jk}=
\sum_{q\in\Lambda}u_{il}^*(u_{q1}u_{q1}^*\cdots
u_{qn}u_{qn}^* )u_{jk}=u_{il}^*u_{i1}u_{i1}^*\cdots
u_{in}u_{in}^* u_{jk}=0
\]
 by Lemma~\ref{lem:4.6}(a). 

We next show that $pu_{ik}\top pu_{il}$ for $k\ne l$. We have
\begin{eqnarray*} \lefteqn{
pu_{ik}(pu_{ik})^*pu_{il}+pu_{il}(pu_{ik})^*pu_{ik}=}&\\ 
&&pu_{ik}u_{ik}^*pu_{il}+pu_{il}u_{ik}^*pu_{ik}\\
&=&(u_{i1}u_{i1}^*\cdots u_{in}u_{in}^*)u_{ik}u_{ik}^*pu_{il}\\
&+&(u_{i1}u_{i1}^*\cdots u_{in}u_{in}^*)u_{il}u_{ik}^*pu_{ik}\\
&=&(u_{i1}u_{i1}^*\cdots u_{in}u_{in}^*)u_{il}+0
\mbox{ (since $u_{ik}^*u_{il}u_{ik}^*=0$)}\\ 
&=&pu_{il}\mbox{ as required }. 
\end{eqnarray*} 

We next show that for $i\ne j$,
$pu_{jk}\top pu_{ik}$. To this end , we
shall show that $pu_{jk}(pu_{jk})^*pu_{ik}
=0$ and $pu_{ik}(pu_{jk})^*pu_{jk}=pu_{ik}$. In the first place,
\[
pu_{jk}(pu_{jk})^*pu_{ik}=(pu_{jk}u_{jk}^*)(pu_{ik})=
(u_{j1}u_{j1}^*\cdots u_{jn}u_{jn}^*)(u_{i1}u_{i1}^*\cdots u_{in}u_{in}^*)=0.
\] 
In the second place,
\begin{eqnarray*}
pu_{ik}(pu_{jk})^*pu_{jk}&=&pu_{ik}u_{jk}^*(pu_{jk})\\
&=&pu_{ik}u_{jk}^*(\prod_{1\le l\le n,l\ne k}u_{jl}u_{jl}^*)u_{jk}\\
&=&pu_{ik}(u_{jk}^*-u_{j1}^*u_{j1}u_{jk}^*)u_{j2}u_{j2}^*\cdots 
u_{jn}u_{jn}^*u_{jk}\\
&=&pu_{ik}u_{jk}^*(u_{j2}u_{j2}^*\cdots  
u_{jn}u_{jn}^*)u_{ik}\\
&&\ldots\\ 
&=&pu_{ik}u_{jk}^*u_{jk}=p(u_{ik}-u_{jk}u_{jk}^*u_{ik})=pu_{ik}. 
\end{eqnarray*}

Finally we shall show that $\tpc{pu_{jk}}{pu_{jl}}{pu_{il}}=pu_{ik}/2$ for 
$j\ne i$ and $l\ne k$. It suffices to prove
$pu_{jk}u_{jl}^*pu_{il}=0$ and $pu_{il}u_{jl}^*pu_{jk}=pu_{ik}$.

On the one hand,
\begin{eqnarray*}
pu_{jk}u_{jl}^*pu_{il}&=&\sum_{m\in\Lambda}u_{m1}u_{m1}^*\cdots u_{mn}u_{mn}^*
u_{jk}u_{jl}^*pu_{il}\\
&=&u_{j1}u_{j1}^*\cdots u_{jn}u_{jn}^*u_{jk}u_{jl}^*pu_{il}=0,
\end{eqnarray*}
since $u_{jl}^*u_{jk}u_{jl}^*=0$.

On the other hand,
\begin{eqnarray*}
pu_{il}u_{jl}^*pu_{jk}&=&\sum_{m\in\Lambda}pu_{il}u_{jl}^*u_{m1}u_{m1}^*
\cdots u_{mn}u_{mn}^*u_{jk}\\
&=&pu_{il}u_{jl}^*[u_{j1}u_{j1}^*\cdots u_{jn}u_{jn}^*]u_{jk}\\
&&\mbox{(where }
u_{jl}u_{jl}^*\mbox{ and }u_{jk}u_{jk}^*\mbox{ are not present
in the }[\ \cdot\ ])\\
&=&pu_{il}[u_{jl}^*-u_{j1}^{\ast}u_{j1}u^{\ast}_{jl}]u_{j2}u_{j2}^*\cdots 
u_{jn}u_{jn}^*u_{jk}\\
&=&pu_{il}u_{jl}^*u_{j2}u_{j2}^*\cdots u_{jn}u_{jn}^*u_{jk}\\
&=&\cdots\\
&=&pu_{il}u_{jl}^*u_{jk}\\
&=&p(u_{ik}-u_{jk}u_{jl}^*u_{il})\\
&=&pu_{ik}-pu_{jk}u_{jl}^*u_{il}\\
&=&pu_{ik}-\sum_{m\in\Lambda}u_{m1}u_{m1}^*\cdots u_{mn}u_{mn}^*
u_{jk}u_{jl}^*u_{il}\\
&=&pu_{ik}-u_{j1}u^{\ast}_{j1} \cdots u_{jn}u^{\ast}_{jn}u_{jk}u^{\ast}_{jl}u_{il}\\
&=&pu_{ik}\mbox{ since }u_{jl}^*u_{jk}u_{jl}^*=0. 
\end{eqnarray*}

It now follows that the map $y\mapsto py$ is a triple isomorphism,
and hence an isometry, from
the norm closure $U$ of 
$\mbox{sp}_{\CC}\, \{u_{ij}\}$ onto
the norm closure $V$ of 
$\mbox{sp}_{\CC}\, \{pu_{ij}\}$.

We claim that the map $y\mapsto py$ is an isometry of the w*-closure 
$\overline{U}$ of
$U$ onto the w*-closure $\overline{V}$ of $V$,
 and is thus a complete contraction as well.

 First we show that if $py=0$ and 
$y\in Y$, then $y=0$, from which it follows 
 that the map $y\mapsto py$ is a w*-homeomorphism when 
restricted to the unit ball of $U$. 
Then by \cite[(3.1)]{Horn87}, $y\mapsto py$ extends to an isometry
of $\overline{U}$ onto $\overline{V}$, which is w*-continuous by the
uniqueness of the preduals.  This w*-extension must agree with 
$y\mapsto py$ on $Y$, which proves the claim. 

To prove the above statement suppose $py=0$ for some $y\in Y$.
Then $L_iy=0$ for each $i\in\Lambda$. We may
write $y=\sum\lambda_{ij}u_{ij}$ where the sum converges in the w*-topology
and $\lambda_{ij}=l_{ij}yr_{ij}$, where $l_{ij}=u_{ij}u_{ij}^*$, 
(resp. $r_{ij}=u_{ij}^*u_{ij}$)
 is the left (resp. right) support 
of $u_{ij}$. Since $L_iu_{ij}R_j=L_iu_{ij}\ne 0$, 
$0=L_iyR_j=L_il_{ij}yr_{ij}R_j=\lambda_{ij}L_iu_{ij}R_j$ and
so $\lambda_{ij}=0$ and $y=0$.

We can similarly 
show that $\{(1-p)u_{ij}\}$ is a rectangular grid and that hence,
as above, the map $Y\ni y\mapsto (1-p)y\in (1-p)Y$ 
is an isometry and complete contraction. For example, to prove
that $(1-p)u_{jk}\top(1-p)u_{ik}$, it suffices to show that
\[
(1-p)u_{jk}[(1-p)u_{ik}]^*(1-p)u_{ik}=0
\]
and
\[
(1-p)u_{ik}[(1-p)u_{ik}]^*(1-p)u_{jk}=(1-p)u_{jk}.
\]  

For the first statement,
\begin{eqnarray*}\lefteqn{
(1-p)u_{jk}u_{ik}^*(1-p)u_{ik}=
(1-p)u_{jk}u_{ik}^*u_{ik}
-(1-p)u_{jk}u_{ik}^*pu_{ik}
}&\\
&=&(1-p)u_{jk}u_{ik}^*u_{ik}-(1-p)u_{jk}u_{ik}^*(u_{i1}u_{i1}^*\cdots
u_{in}u_{in}^*)u_{ik}\\
&=&(1-p)u_{jk}u_{ik}^*u_{ik}-(1-p)u_{jk}u_{ik}^*u_{ik}=0,
\end{eqnarray*}
since in the second term,
 for $l\ne k$, $u_{jk}u_{ik}^*u_{il}u_{il}^*=u_{jk}(u_{ik}^*
-u_{il}^*u_{il}u_{ik}^*)=u_{jk}u_{ik}^*$. 

For the second statement,
\begin{eqnarray*} 
(1-p)u_{ik}u_{ik}^*(1-p)u_{jk}&=&(1-p)u_{ik}u_{ik}^*u_{jk}-(1-p)u_{ik}
u_{ik}^*pu_{jk}\\
&=&(1-p)u_{ik}u_{ik}^*u_{jk}+0\\
&=&(1-p)(u_{jk}-u_{jk}u_{ik}^*u_{ik})=(1-p)u_{jk}.
 \end{eqnarray*}

We omit the entirely analogous calculations showing that the other
grid properties hold. 

As above, the fact that $y\mapsto (1-p)y$ is a 
complete semi-isometry follows from the fact that it is one-to-one on $Y$.
To see that it is one-to-one on $Y$, suppose $(1-p)y=0$ for some $y\in Y$.
Writing $y=\sum\lambda_{ij}u_{ij}$ leads to $\sum_{i,j}\lambda_{ij}(1-L_i)u_{i,j}=0$.
If there were indices $(i,j)$ such that $(1-L_i)u_{i,j}=0$, then since $(1-L_i)u_{i,j}=(1-p)u_{ij}$,
we would have $u_{kj}^*u_{ij}=u_{kj}^*pu_{ij}=u_{kj}^*L_iu_{ij}=0$ for some $k$, violating
(\ref{eq:15bis}).

Finally we show that $pY\perp (1-p)Y$. It suffices to show that basis
elements are orthogonal. First, 
\[
pu_{ij}[(1-p)u_{ij}]^*=pu_{ij}u_{ij}^*(1-p)
=(\prod_k u_{ik}u_{ik}^*)u_{ij}u_{ij}^*(1-p)=p(1-p)=0. 
\] 

Next, if $j\ne k$, $pu_{ij}u_{ik}^*(1-p)=(\prod_l u_{il}u_{il}^*)
u_{ij}u_{ik}^*(1-p)=0$, since $u_{il}u_{ij}^*u_{il}=0$.

Finally, if $k\ne i$,
\begin{eqnarray*}
pu_{ij}u_{kj}^*(1-p)&=&pu_{ij}u_{kj}^*-pu_{ij}u_{kj}^*\prod_{l=1,l\ne j}^n
u_{kl}u_{kl}^*\\
&=&pu_{ij}u_{kj}^*-pu_{ij}[u_{kj}-u_{k1}u_{k1}^*u_{kj}]^*\prod_{l=2,l\ne j}^n
u_{kl}u_{kl}^*\\
&=&pu_{ij}u_{kj}^*-pu_{ij}u_{kj}^*\prod_{l=2,l\ne j}^n
u_{kl}u_{kl}^*\\
&&\ldots\\ 
&=&pu_{ij}u_{kj}^*-pu_{ij}u_{kj}^*=0. 
\end{eqnarray*}

Clearly $[(1-p)y]^*pz=0$ for all $y,z\in Y$, finishing the proof. \qed

\begin{prop}\label{prop:5.7}
Suppose that $Y$ is triple isomorphic to $B(H,K)$ and is of rank
at least 2, and that {\rm (\ref{eq:466})} and {\rm (\ref{eq:15bis})} hold. 
Then
$Y$ is completely semi-isometric to $B(H,K)$ and completely isometric to
\[
\mbox{Diag}\, (B(H,K),B(K,H)).
\]
\end{prop}
\pf\
For $k\ne j$,
\begin{eqnarray*}
pu_{ik}(pu_{ij})^*&=&pu_{ik}u_{ij}^*p\\
&=&pu_{ik}u_{ij}^*u_{i1}u_{i1}^*\cdots u_{in}u_{in}^*=0.
\end{eqnarray*}
so Lemma~\ref{lem:4.5} applies to show that $pY$ is completely isometric
to $B(H,K)$.  By Lemma~\ref{lem:5.5}, $Y$ is completely semi-isometric to 
$B(H,K)$.

Similarly,  Lemma~\ref{lem:4.5} applies to show that 
$(1-p)Y$ is completely isometric to $B(K,H)$. 
Indeed,
\begin{eqnarray*}
[(1-p)u_{ik}][(1-p)u_{jk}]^*&=&(1-p)u_{ik}u_{jk}^*(1-p)\\
&=&(1-\prod_l u_{il}u_{il}^*)u_{ik}u_{jk}^*(1-p)\\
&=&u_{ik}u_{jk}^*(1-p)-
(\prod_l u_{il}u_{il}^*)u_{ik}u_{jk}^*(1-p)\\
&=&u_{ik}u_{jk}^*(1-p)-
u_{ik}u_{jk}^*(1-p)=0,
\end{eqnarray*}
since for $l\ne k$, $u_{il}u_{il}^*u_{ik}u_{jk}^*=(u_{ik}-u_{ik}
u_{il}^*u_{il})u_{jk}^*=u_{ik}u_{jk}^*$. 

As in the proof of 
Lemma~\ref{lem:3.2}, $pY\oplus^{\ell^\infty}(1-p)Y$ is completely
isometric to $B(H,K)\oplus^{\ell^\infty} B(K,H)$. \qed

\medskip

This completes the proof of Proposition~\ref{prop:1} and Theorem~\ref{thm:3}
in the case
that $Y$ is of rank 2 or more and triple isomorphic to $B(H,K)$.

\subsection{The case of rank 1. Preliminary cases}\label{5.3}

Assume now that $Y$ is finite dimensional and rank 1.
Let us denote a finite rectangular grid of rank 1 for $Y$ by $\{u_1,\ldots,u_n\}$.
For the record, let us note that this means precisely that
$u_i$ is a non-zero partial isometry: $\tp{u_i}{u_i}{u_i}=u_i\ne 0$;
$u_i$ is minimal:
\begin{equation}\label{eq:min}
\tp{u_i}{u_j}{u_i}=0\mbox{ for }i\ne j;
\end{equation}
 and that $u_i$ is 
colinear with $ u_k$: $\tp{u_i}{u_i}{u_k}=u_k/2$ for $i\ne k$. By the grid properties and the identity $\| yy^{\ast}y
\|=\|y\|^{3}$, $Y$ is isometric to a Hilbert space with orthonormal basis $\{ u_{j} \}$ (see \cite[p.306]{DanFri87}).

We shall denote, for $J=\{j_1,\ldots,j_i\}\subset\{1,2,\ldots,n\}$,
$u_{j_1}^*u_{j_1}u_{j_2}^*u_{j_2}\cdots u_{j_i}^*u_{j_i}$ by $(u^*u)_J$. By
commutativity of the projections $u_k^*u_k$ we may and shall assume that
$1\le j_1<\cdots<j_i\le n$. Similarly $(uu^*)_J$ will denote
$u_{j_1}u_{j_1}^*u_{j_2}u_{j_2}^*\cdots u_{j_i}u_{j_i}^*$.

\begin{lem}\label{lem:4.7}
If $(uu^*)_J=0$ for some $J$ with $|J|=i$, then 
$(uu^*)_J=0$ for all $J$ with $|J|=i$.
If $(u^*u)_J=0$ for some $J$ with $|J|=i$, then 
$(u^*u)_J=0$ for all $J$ with $|J|=i$.
\end{lem}
\pf\
Suppose that $(uu^*)_J=0$ for some $J$ with $|J|=i$.
 Then for $s\in J$ and $k\not\in J$, 
\begin{eqnarray*}
(uu^*)_{(J-\{s\})\cup\{k\}}&=&(uu^*)_{J-\{s\}}u_ku_k^*\\
&=&(uu^*)_{J-\{s\}}(u_ku_s^*u_s+u_su_s^*u_k)u_k^*\\
&=&(uu^*)_{J-\{s\}}u_ku_s^*u_su_k^*\\
&=&u_ku_s^*(uu^*)_{J-\{s\}}u_su_k^*
\mbox{ (by Lemma~\ref{lem:4.6}) }\\
&=&u_ku_s^*u_su_s^*(uu^*)_{J-\{s\}}u_su_k^*\\
&=&u_ku_s^*(uu^*)_Ju_su_k^*=0.
\end{eqnarray*}

The proof of the second statement is similar.\qed

\medskip

Lemma~\ref{lem:4.7} makes it possible to define $i_R$ to be 
the largest $i$ such that $(uu^*)_J\ne 0$ for any $J$ with $|J|=i$ and
$i_L$ to be the largest $i$ such that 
$(u^*u)_J\ne 0$ for any $J$ with $|J|=i$. The numbers $i_R$ and $i_L$ are
indices which depend on how a \jcst\ sits in its ternary envelope.
We use the numbers $i_R$ and $i_L$ to define projections $p_R=\sum
_{|J|=i_R}(uu^*)_J$ and $p_L=\sum_{|J|=i_L}(u^*u)_J$.

\begin{lem}\label{lem:5.8} Each of
the maps $y\mapsto p_Ry,\ y\mapsto yp_L,\ y\mapsto (1-p_R)y,\ y\mapsto 
y(1-p_L)$
 are completely contractive triple
isomorphisms of $Y$ into~ $A$. Also $p_RY\perp (1-p_R)Y$ and $Yp_L\perp
Y(1-p_L)$.
\end{lem}
\pf\
To prove the first statement, it 
suffices to show that each of these maps takes the rectangular rank 1 grid
$\{u_k\}_{k=1}^n$ into a rectangular grid of rank~ 1.

We carry out the proof for $y\mapsto p_Ry$ and $y\mapsto (1-p_R)y$,
the proofs for the other maps being identical.
For notation's sake, we let $p=p_R$ and $w_j=pu_j$. 

If $w_k=0$, then 
\[
0=w_kw_k^*=(pu_k)(pu_k)^*=pu_ku_k^*p=
\sum_{|J|=i_R,k\in J}(uu^*)_J. 
\]
A sum of orthogonal projections cannot be
zero unless each one is. Thus $(uu^*)_J=0$ for any $J$ containing $k$ with
$|J|=i_R$, which is a 
contradiction. Hence, $w_k\ne 0$.

Next, for $i\ne j$,
\begin{eqnarray*}
w_iw_j^*w_i&=&pu_i(pu_j)^*pu_i=pu_iu_j^*pu_i\\
&=&pu_iu_j^*
\left(\sum_{|J|=i_R}(uu^*)_J\right)u_i=0,
\end{eqnarray*}
since if $i\in J$, $u_iu_j^*(uu^*)_J=0$ and if $i\not\in J$, then
$(uu^*)_Ju_i=0$.

Similarly,
\[    w_iw_i^*w_i=
pu_i(pu_i)^*pu_i=pu_i^*u_ipu_i=\left(\sum_{|J|=i_R,i\in J}(uu^*)_J\right)
u_i
\]
and $w_i=
pu_i=\sum_{|J|=i_R}(uu^*)_Ju_i=\sum_{|J|=i_R,i\in J}(uu^*)_Ju_i$, so that
$w_iw_i^*w_i=w_i$.

Now we shall show that $w_i$ and $w_k$ are colinear. 
It suffices to show that for $i\ne k$,
\[
pu_i(pu_i)^*pu_k+pu_k(pu_i)^*pu_i=pu_k,
\]
equivalently (by using Lemma~\ref{lem:4.6} on the middle term),
\[
pu_iu_i^*pu_ku_k^*+pu_ku_i^*pu_iu_k^*=pu_ku_k^*.
\]
As noted above, 
\begin{equation}\label{eq:16prime}
pu_ku_k^*=pu_ku_k^*p=
\sum_{|J|=i_R,k\in J}(uu^*)_J.
\end{equation}
 On the other
hand we have
\begin{equation}\label{eq:11bis}
pu_iu_i^*pu_ku_k^*=\sum_{|J|=i_R,i,k\in J}(uu^*)_J.
\end{equation}
and
\begin{equation}\label{eq:661}
pu_ku_i^*pu_iu_k^*=
\left(\sum_{|J|=i_R,i\not\in J}(uu^*)_J\right)u_ku_i^*\left(
\sum_{|J|=i_R,k\not\in J}(uu^*)_J\right)u_iu_k^*.
\end{equation}

It remains to show that the right side of (\ref{eq:661}), call it
$A$, when added to the right side of (\ref{eq:11bis}), equals the
right side of (\ref{eq:16prime}). 

We have 
\begin{eqnarray*}
A&=&\left(\sum_{|J_1|=i_R,i\not\in J_1}(uu^*)_{J_1}\right)
u_ku_i^*
\left(\sum_{|J_2|=i_R,k\not\in J_2}(uu^*)_{J_2}\right)u_iu_k^*\\
&=&
\sum_{|J_1|=|J_2|=i_R,i\not\in J_1,k\not\in J_2}
(uu^*)_{J_1}u_ku_i^*(uu^*)_{J_2}u_iu_k^*\\
\end{eqnarray*}

Now each term in this sum for which $k\not\in J_1$ is zero, as is each
term for which $i\not\in J_2$.  On the other hand, if
 $i\in J_2$ and $k\in J_1$, then
by Lemma~\ref{lem:4.6},
\[
(uu^*)_{J_1}u_ku_i^*(uu^*)_{J_2}u_iu_k^*=
(uu^*)_{J_1}(uu^*)_{J_2-\{i\}}u_ku_i^*u_iu_k^*,
\]
which is zero unless $J_1=J_2\cup\{k\}$, in which case it
equals $(uu^*)_{J_1}u_ku_i^*u_iu_k$, 
where $k\in J_1$ and $i\not\in J_1$.

Conversely, if $k\in J$ and $i\not\in J$, then 
\[
(uu^*)_Ju_ku_i^*u_iu_k^*=(uu^*)_{J_1}u_ku_i^*(uu^*)_{J_2}u_iu_k^*
\]
where $J_1=J,\ J_2=(J-\{k\})\cup\{i\}$, $|J_1|=|J_2|=i_R$, $i\not\in J_1,k\not\in J_2$.
  Therefore
\begin{eqnarray*}
A&=&\sum_{|J|=i_R,k\in J,i\not\in J}(uu^*)_Ju_ku_i^*u_iu_k^*\\
&=&\sum_{|J|=i_R,k\in J,i\not\in J}(uu^*)_J(u_k-u_iu_i^*u_k)u_k^*\\
&=&\sum_{|J|=i_R,k\in J,i\not\in J}(uu^*)_Ju_ku_k^*\\
&=&\sum_{|J|=i_R,k\in J,i\not\in J}(uu^*)_J,
\end{eqnarray*}
as required. This proves that $\{p_Ru_k\}_{k=1}^n$ is a rectangular rank
1 grid.

Let us now prove that $(1-p)u_i\top (1-p)u_k$, that is,
\[
(1-p)u_i[(1-p)u_i]^*(1-p)u_k\ +\ (1-p)u_k[(1-p)u_i]^*(1-p)u_i=(1-p)u_k.
\] 
As  before, it suffices to prove
\begin{equation}\label{eq:8171}
(1-p)u_iu_i^*(1-p)u_ku_k^*+(1-p)u_ku_i^*(1-p)u_iu_k^*=(1-p)u_ku_k^*.
\end{equation}
For the first term on the left side of (\ref{eq:8171}),
\begin{eqnarray}\label{eq:8172}\nonumber
\lefteqn{(1-p)u_iu_i^*(1-p)u_ku_k^*=}&\\
&=&(1-p)u_iu_i^*u_ku_k^*- 
(1-p)u_iu_i^*pu_ku_k^*\\
\nonumber
&=&(1-p)u_iu_i^*u_ku_k^*\mbox{ (since $u_iu_i^*$ commutes with $p$)}
\end{eqnarray}
For the second term on the left side of (\ref{eq:8171}),
\begin{eqnarray}\nonumber\label{eq:8173}
\lefteqn{(1-p)u_ku_i^*(1-p)u_iu_k^*=}&\\\nonumber
&=&   (1-p)u_ku_i^*u_iu_k^*-(1-p)u_ku_i^*pu_iu_k^*\\\nonumber
&=&(1-p)u_ku_i^*u_iu_k^*-(1-p)u_ku_i^*\left(\sum_{|J|=i_R}(uu^*)_J\right)
u_iu_k^*\\\nonumber
&=&(1-p)u_ku_i^*u_iu_k^*-(1-p)u_ku_i^*\left(\sum_{|J|=i_R,i\in J,k\not\in J}
(uu^*)_J\right)u_iu_k^*\\
&=&(1-p)u_ku_i^*u_iu_k^*-(1-p)u_ku_i^*\left(\sum_{|J|=i_R-1,i\not\in J,k\not\in J}
(uu^*)_J\right)u_iu_k^*\\\nonumber
&=&(1-p)u_ku_i^*u_iu_k^*-(1-p)\left(\sum_{|J|=i_R-1,i\not\in J,k\not\in J}
(uu^*)_J\right)u_ku_i^*u_iu_k^*\\\nonumber
&=&(1-p)u_ku_i^*u_iu_k^*-(1-p)\left(\sum_{|J|=i_R,i\not\in J,k\in J}
(uu^*)_J\right)u_ku_i^*u_iu_k^*\\\nonumber
&=&(1-p)u_ku_i^*u_iu_k^*-(1-p)p
u_ku_i^*u_iu_k^*\\\nonumber
&=&(1-p)u_ku_i^*u_iu_k^*. 
\end{eqnarray} 

By (\ref{eq:8172}) and (\ref{eq:8173}), the left side of (\ref{eq:8171})
is equal to 
\[
(1-p)u_ku_k^*u_iu_i^*+(1-p)u_ku_i^*u_iu_k^*=
(1-p)u_k(u_k^*u_iu_i+u_i^*u_iu_k^*)=(1-p)u_ku_k^*,
\]
as required.

We omit the analogous proof that $(1-p)u_i$ is a minimal partial isometry.

Finally we show that $pY\perp (1-p)Y$. It suffices to show that basis
elements are orthogonal, that is $pu_i[(1-p)u_j]^*=0$ for all $i,j$. First,
if $i\ne j$, then
\begin{eqnarray*}
pu_i[(1-p)u_j]^*&=&pu_iu_j^*(1-p)\\
&=&\left(\sum_{i,j\not\in J,|J|=i_R-1}(uu^*)_J\right)u_iu_j^*(1-p)\\
&=&u_iu_j^*\left(\sum_{i,j\not\in J,|J|=i_R-1}(uu^*)_J\right)(1-p)\\
&=&u_iu_j^*\left(\sum_{i\not\in J,|J|=i_R}(uu^*)_J\right)(1-p)\\
&=&u_iu_j^*p(1-p)=0.
\end{eqnarray*}

Next, $pu_i[(1-p)u_i]^*=pu_iu_i^*-pu_iu_i^*p=pu_iu_i^*-pu_iu_i^*=0$.

Clearly, $(pz)^*(1-p)y=0$ for all $y,z\in Y$.\qed

\medskip

The next proposition proves Theorem~\ref{thm:1}(c) (See 
Remark~\ref{rem:6.05} for the definition of the spaces $H_n^k$). 

\begin{prop}\label{prop:5.10bis}
If either of $i_R$ or $i_L$ is equal to 1 or $n$, then $Y$ is completely
semi-isometric to $R_n$ or to $C_n$.
\end{prop}
\pf\
If $i_R=1$, then $u_1^*u_2=0$ so by Corollary~\ref{cor:1}(b), $Y$ is
completely isometric to $B(\CC^n,\CC)$.  If $i_L=1$, then $u_1u_2^*=0$ so
by Corollary~\ref{cor:1}(a), $Y$ is completely isometric to 
$B(\CC,\CC^n)$.

If $i_R=n$, then with $p=p_R$,
\[
pu_1(pu_2)^*=pu_1u_2^*p
=(uu^*)_{\{1,2,\ldots,n\}}
u_1u_2^*p=0
\]
since $u_2^*u_1u_2^*=0$, so by Corollary~\ref{cor:1}(a), $pY$
is completely isometric to $B(\CC,\CC^n)$ and by Lemma~\ref{lem:5.8},
 $Y$ is completely
semi-isometric to $B(\CC,\CC^n)$. Similarly, if $i_L=n$, then
$Y$ is completely semi-isometric to $B(\CC^n,\CC)$.
\qed


\medskip

In preparation for the next two sections, let us
 consider the remaining case where $1<i_R,i_L<n$. 

\begin{lem}\label{lem:5.10} In general,
$i_R+i_L\ge n+1$.
Let $p=p_R$ and $w_j
=pu_j$. Let $i_L^\prime$ and $i_R^\prime$ denote the corresponding indices
for the grid $\{w_1,\ldots,w_n\}$.
Then 
$i_L^\prime+i_R^\prime =n+1$.
\end{lem}
\pf\
Note first that if $i_L<n$,
\begin{eqnarray}\label{eq:719bis}\nonumber 
(u^*u)_{\{1,2,\ldots,i_L\}}&=&
(u^*u)_{\{1,2,\ldots,i_L-1\}}u_{i_L}^*u_{i_L}\\\nonumber
&=&(u^*u)_{\{1,2,\ldots,i_L-1\}}u_{i_L}^*
(u_{i_L}u_{i_L+1}^*u_{i_L+1}+
u_{i_L+1}u_{i_L+1}^*u_{i_L})\\\nonumber
&=&0+(u^*u)_{\{1,2,\ldots,i_L-1\}}u_{i_L}^*u_{i_L+1}u_{i_L+1}^*
u_{i_L}\\\nonumber
&=&(u^*u)_{\{1,2,\ldots,i_L-1\}}u_{i_L}^*u_{i_L+1}u_{i_L+1}^*
(u_{i_L}u_{i_L+2}^*u_{i_L+2}+u_{i_L+2}u_{i_L+2}^*u_{i_L})\\
&=&(u^*u)_{\{1,2,\ldots,i_L-1\}}u_{i_L}^*
(u_{i_L+1}u_{i_L+1}^*
u_{i_L+2}u_{i_L+2}^*)u_{i_L})\\\nonumber
&=&\ldots\\\nonumber
&=&(u^*u)_{\{1,2,\ldots,i_L-1\}}u_{i_L}^*(uu^*)_
{\{i_L+1,i_L+2,\ldots,n\}}u_{i_L}.\\ \nonumber
&=&(u^*u)_{\{1,2,\ldots,i_L-1\}}u_{i_L}^*(uu^*)_
{\{i_L,i_L+1,i_L+2,\ldots,n\}}u_{i_L}.
\end{eqnarray}
If $n-i_L+1>i_R$, then 
$|\{i_L,i_L+1,\ldots,n\}|>i_R$, so 
$(uu^*)_{\{1,2,\ldots,i_L \}}=0$, which is impossible.
Hence $i_L+i_R\ge n+1$, proving the first statement.

It is easy to see that $(w^*w)_J=(u^*u)_J$ and therefore that $i_L^\prime
=i_L$. Indeed, for any $r\ge 1$,
\begin{eqnarray*}
(w^*w)_{\{1,2,\ldots,r\}}&=&pu_1^*u_1pu_2^*u_2\cdots pu_r^*u_rp\\
&=&\sum_{|J|=i_L,\{1,\ldots,r\}\subset J}(u^*u)_J.
\end{eqnarray*}

Moreover, for any $r\ge 1$,
\begin{eqnarray}\label{eq:622}\nonumber
(ww^*)_{\{1,\ldots,r\}}&=&u_1pu_1^*u_2pu_2^*\cdots u_rpu_r^*\\
&=&
\sum[u_1(u^*u)_{J_1}u_1^*]\cdots
[u_r(u^*u)_{J_r}u_r^*],
\end{eqnarray}
where the sum can be taken over all $|J_k|=i_L$ with $k\in J_k$ and
$(\{1,2,\ldots,r\}- \{k\})\cap J_k=\emptyset$. 
Indeed, if $k\not\in J_k$, then $u_k(u^*u)_{J_k}u_k=0$; and if there is
a $j\in (\{1,2,\ldots,r\}- \{k\})\cap J_k$, then by commutativity
of the factors in the terms of (\ref{eq:622}), that term would vanish
by (\ref{eq:min}).

Thus if $r\le i_R^\prime$, we have $i_L=|J_k|\le n-(r-1)$ and in
particular, $i_L\le n-(i_R^\prime-1)$, that is, $i_L+i_R^\prime\le
n+1$.
Since $i_L=i_L^\prime$ and $i_R^\prime +i_L^\prime \ge n+1$, we conclude
that $i_L^\prime+i_R^\prime =n+1$.\qed

\section{The Hilbertian operator spaces $H_n^k$}

In this section we shall begin by assuming that $Y$  is a
\jwst\ of rank 1 and 
finite dimension $n$ given
by a rectangular rank 1 grid $\{u_1,\ldots,u_n\}$ such
that $i_R+i_L=n+1$.  If $i_L=1$ or if $i_R=1$, then $Y$ is completely
isometric to the type~ 1 Cartan factors $R_n$ or $C_n$ by 
Corollary~\ref{cor:1}. Otherwise, 
we shall show in section 7 that $Y$ is completely isometric to
a space $H_n^{i_R}$ 
which is a subtriple of a Cartan factor of type~ 1, proving
Proposition~\ref{prop:1} in this case. This will be
achieved by constructing, from the given grid $\{u_j\}$,
a rectangular grid $\{u_{IJ}\}$ whose linear span is a 
ternary algebra containing $Y$ and which is ternary isomorphic to a 
Cartan factor of type 1, namely the $n\choose i_L$ by $n\choose i_R$ complex
matrices. 

After this, in section 7 we shall prove Proposition~\ref{prop:1} in case
$Y$ is infinite dimensional and of rank 1.

Here is the construction. 
We define some elements which are indexed by an 
arbitrary pair of subsets $I,J$ of $\{1,\ldots,n\}$
satisfying 
\begin{equation}\label{eq:12}
|I|=i_R-1,|J|=i_L-1.
\end{equation} 
Note that the number of possible sets $I$ is $n\choose{i_R-1}$
($={n\choose {i_L}}$) and the number of such $J$ is
$n\choose{i_L-1}$
($={n\choose {i_R}}$). Moreover, if $|I\cap J|=s\ge 0$, 
then $|(I\cup J)^c|=s+1$. 
 Hence we may write
\[
I=\{i_1,\ldots,i_k,d_1,\ldots,d_s\}\quad , \quad J=\{j_1,\ldots,j_l,d_1,\ldots,d_s\}
\]
where $I\cap J=\{d_1,\ldots,d_s\}$.
Let us write $ (I\cup J)^c=\{c_1,\ldots,c_{s+1}\}$, and
let us agree (for the moment) that the elements are ordered
as follows: $c_1<c_2<\cdots <c_{s+1}$ and $d_1<d_2<\cdots<d_s$. 

\begin{defn}\label{def:6.0}
With the above notation, we define
\begin{equation}\label{eq:(I,J)}
u_{IJ}=u_{I,J}=(uu^*)_{I- J}u_{c_1}u_{d_1}^*u_{c_2}u_{d_2}^*\cdots u_{c_s}u_{d_s}^*u_{c_{s+1}}(u^*u)_{J- I}.
\end{equation}
\end{defn}

\begin{rem}\label{rem:6.05}
We   are going to show (c.f. Proposition \ref{prop:4.14} and \ref{prop:4.17}) that there is a choice of signs $\epsilon(I,J)
=\pm 1$ such that the family $\{\epsilon(I,J)u_{I,J}\}$
forms a rectangular grid which is closed under the ternary 
product $(a,b,c)\mapsto ab^*c$, so that its linear span
is ternary isomorphic and
therefore completely isometric 
to a concrete Cartan factor of type 1. By restriction, from
{\rm (\ref{eq:sum})} below, 
 $Y$ will be completely isometric to its image,
which we shall denote by $H_n^{i_R}$. We will then show that all $H_n^{k}$ are actually rank 1 triples (and thus
Hilbertian) and satisfy $k= i_{R}, i_{R}+i_{L}=n+1$, thus proving the existence of the Hilbert spaces discussed in
this section (see the paragraph preceding Example 1 in section 7).
\end{rem}
\begin{prop}\label{prop:4.14}
Let $u_{I,J}$ be defined by {\rm (\ref{eq:(I,J)})}. Then
\begin{description}
\item[(a)] $u_{I,J}$ is a minimal partial isometry, that is
\[
u_{I,J}[u_{I,J}]^*u_{I,J}=u_{I,J}\mbox{ and }
u_{I,J}[u_{I',J'}]^*u_{I,J}=0\mbox{ for all }(I,J)\ne (I',J').
\] 
\item[(b)] orthogonality: $u_{I,J}\perp u_{I',J'}$ if 
$I\ne I'$ and $J\ne J'$.
\item[(c)] colinearity: $u_{I,J}\top u_{I',J'}$ if either $I=I'$ or
$J=J'$ {\rm (}but not both{\rm )}.
\item[(d)] associative orthogonality: 
\[
u_{I,J}[u_{I',J'}]^*=0\mbox{ if }I\ne I';\mbox{ and }
[u_{I,J}]^*u_{I',J'}=0\mbox{ if }J\ne J'.
\] 
\item[(e)] ``weak'' quadrangle property: 
$u_{I,J}[u_{I,J'}]^*u_{I',J'}=\pm u_{I',J}$
\end{description}
\end{prop}
\pf\
Throughout this proof, we use the fact that all elements of the grid
$\{u_1,\ldots,u_n\}$ are present in each $u_{IJ}$. To avoid 
cumbersome notation we
will also
often denote an element $u_{c}$, where $c\in (I\cup J)^c$, by
 $c_{ij}$, and similarly for
$u_{d}$. For example, in (\ref{eq:623}) below, $d^1_{ij'}$ denotes
$u^*_{d^1_{ij'}}$ where $d^1_{ij'}\in I\cap J'$, and $c_{ij}^1$
denotes $u_{c_{ij}^1}$, where $c_{ij}^1\in (I\cup J)^c$.
 
\medskip

{\bf Proof of (e):} By definition, 
\begin{eqnarray}\nonumber 
u_{I,J}[u_{I,J'}]^*u_{I',J'}&=&
\left[(uu^*)_{I- J}c_{ij}^1d_{ij}^1\cdots d_{ij}^qc_{ij}^{q+1}
(u^*u)_{J- I}\right]\\\label{eq:623}
&\times&
\left[
(u^*u)_{J'- I}c_{ij'}^{r+1}d_{ij'}^r\cdots d_{ij'}^1c_{ij'}^1
(uu^*)_{
I- J'}
\right]\\\nonumber
&\times&
\left[
(uu^*)_{I'- J'}c_{i'j'}^1d_{i'j'}^1\cdots d_{i'j'}^sc_{i'j'}^{s+1}
(u^*u)_{J'- I'}
\right].
\end{eqnarray}

This quantity remains unchanged if the factors 
\[
(u^*u)_{J- I}(u^*u)_{J'- I}
\mbox{ and }
(uu^*)_{
I- J'}
(uu^*)_{I'- J'}
\]
are removed.   Indeed, since $J- I\subset (I^c\cap {J'}^c)
\cup (I^c\cap J'$ (disjoint union), by using Lemma~\ref{lem:4.6} 
$(u^*u)_{J- I}$ can be absorbed into the $c_{ij'}$'s or into 
$(uu^*)_{J'- I}$.
Similarly, 
$(uu^*)_{I'- J'}$ can be absorbed into the $c_{ij'}$'s or into
$(uu^*)_{I- J'}$.
After this has been done, 
$(u^*u)_{J'- I}$ can be absorbed into the $d_{i'j'}$'s or into
$(u^*u)_{J'- I'}$, and
 $(uu^*)_{I- J'}$ can be absorbed into the $d_{ij}$'s or into
$(uu^*)_{I- J}$,

Thus
\begin{eqnarray}\label{eq:111}\nonumber
\lefteqn{u_{I,J}[u_{I,J'}]^*u_{I',J'}=}\\
&&(uu^*)_{I- J}c_{ij}^1d_{ij}^1\cdots d_{ij}^qc_{ij}^{q+1}
c_{ij'}^{r+1}d_{ij'}^r\cdots d_{ij'}^1c_{ij'}^1
c_{i'j'}^1d_{i'j'}^1\cdots d_{i'j'}^sc_{i'j'}^{s+1}
(u^*u)_{J'- I'}.
\end{eqnarray} 

We claim next that in fact 
\begin{eqnarray}\label{eq:222}\nonumber 
\lefteqn{u_{I,J}[u_{I,J'}]^*u_{I',J'}=}\\
&&(uu^*)_{I'- J}c_{ij}^1d_{ij}^1\cdots d_{ij}^qc_{ij}^{q+1}
c_{ij'}^{r+1}d_{ij'}^r\cdots d_{ij'}^1c_{ij'}^1
c_{i'j'}^1d_{i'j'}^1\cdots d_{i'j'}^sc_{i'j'}^{s+1}
(u^*u)_{J- I'}.
\end{eqnarray}

To get from (\ref{eq:111}) to (\ref{eq:222}) we proceed as follows. 
Consider first an element $x\in I'- J$. Either $x\in I$ or
$x\not\in I$. In the latter case, $x\in (I\cup J)^c$ so that $u_x$
is one of the $c_{ij}$ and so $u_xu_x^*$ can be 
split off from $u_x=u_xu_x^*u_x$
and absorbed (using Lemma~\ref{lem:4.6}) into 
the $(uu^*)_{I- J}$ term. 
In the former case, no
absorption is necessary.  Doing this for every such $x$ allows us
to replace the term
$(uu^*)_{I- J}$ in (\ref{eq:111}) by $(uu^*)_{(I\cup I')
- J}$.

We now have 
\begin{eqnarray}\label{eq:333}\nonumber
\lefteqn{u_{I,J}[u_{I,J'}]^*u_{I',J'}=}\\ &&
(uu^*)_{(I\cup I')- J}c_{ij}^1d_{ij}^1\cdots d_{ij}^qc_{ij}^{q+1}
c_{ij'}^{r+1}d_{ij'}^r\cdots d_{ij'}^1c_{ij'}^1
c_{i'j'}^1d_{i'j'}^1\cdots d_{i'j'}^sc_{i'j'}^{s+1}
(u^*u)_{J'- I'}.
\end{eqnarray} 

Now consider an element $x\in I- J$. Either $x\in J'$ or
$x\not\in J'$.  In the first case, $x\in I\cap J'$ so that $u_x$ is one
of the $d_{ij'}$ and  therefore any such $u_xu_x^*$ can be absorbed from the term $(uu^*) 
_{(I\cup I')- J}$ into a $d_{ij'}$.  On the other hand, if $x\not \in J'$, then
either $x\in I'$, in which case no absorption is necessary, or
$x\not\in I'$ so that $x\in (I'\cup J')^c$ and $u_x$ is one of the
$c_{i'j'}$ and  hence $u_xu_x^*$ can be absorbed.
Doing this for every such $x$ allows us
to replace the term
$(uu^*)_{(I\cup I')- J}$ in (\ref{eq:333}) by $(uu^*)_{ I'
- J}$. 

By an entirely similar two-step argument, we may replace $(u^*u)_{J'-
I'}$ in (\ref{eq:111}) by 
$(u^*u)_{J-
I'}$, which proves (\ref{eq:222}).

\medskip

To complete the proof of (e), we need to show that the right side of 
(\ref{eq:222}) has the form
\[
\pm (uu^*)_{I'- J}c_{i'j}^1d_{i'j}^1
\cdots d_{i'j}^tc_{i'j}^{t+1}(u^*u)_{J- I'}.
\]

To do this we must examine each of the elements
$c_{ij},d_{ij'},c_{i'j'}$ in (\ref{eq:222})
(call them ``outer'' elements as they are
not ``starred'') and $d_{ij},c_{ij'},d_{i'j'}$ 
(call them ``inner'' elements, as they are ``starred'')
and decide whether to leave the element there or absorb it into one of the
end terms $(u^*u)_{J- I'}$ or $(uu^*)_{I'- J}$. This is
achieved in the following lemma. 

\begin{lem} \label{lem:4.15}
Retain the above notation.
\begin{description}
\item[(a)] Each ``outer'' element $c_{ij},d_{ij'},c_{i'j'}$ 
on the right side of {\rm (\ref{eq:222})} 
is either equal to
a $c_{i'j}$ or is equal to a unique other element
on the right side of {\rm (\ref{eq:222})}, together with which it can
be absorbed into one of the terms
$(u^*u)_{J- I'}$ or $(uu^*)_{I'- J}$. Conversely, every
$c_{i'j}$ is equal to one of these ``outer'' elements.
\item[(b)] Similarly, each 
``inner'' element $d_{ij},c_{ij'},d_{i'j'}$ is either equal 
to a $d_{i'j}$ or is equal to a unique other element, together with
which it can
be absorbed into one of the terms
$(u^*u)_{J- I'}$ or $(uu^*)_{I'- J}$.
Conversely, every
$d_{i'j}$ is equal to one of these ``inner'' elements.
\end{description}
\end{lem}
{\bf Proof of Lemma~\ref{lem:4.15}.}
For three mutually colinear partial isometries $u,v,w$, 
the term ``flipping'' in this proof
refers to the
fact that $uv^*w=-wv^*u$.

 Let $c_{ij}\in (I\cup J)^c$.  Either $c_{ij}\in I'$ or $c_{ij}\not\in I'$.
In the first case $c_{ij}$ is a $c_{i'j}$ and no absorption is necessary.
In the second case, either $c_{ij}\in J'$ or $c_{ij}\not\in J'$. If the
former, $c_{ij}\in I'\cap J'$ so that $c_{ij}$ is equal to a $d_{i'j'}$
with which it can be paired by ``flipping'' and $c_{ij}c_{ij}^*$ can
be absorbed into $(uu^*)_{I'- J}$ by Lemma~\ref{lem:4.6}.
In the latter, $c_{ij}\in (I\cup J')^c$ so that  $c_{ij}$ is equal to a
$c_{ij'}$ with which it can be paired and absorbed as above by repeated
use of Lemma~\ref{lem:4.6}.

 Let $d_{ij'}\in I\cap J'$.  Either $d_{ij'}\not \in I'$ or 
$d_{ij'}\in I'$.
In the second case $d_{ij'}\in I'$ so that $d_{ij'}\in I'\cap J'$
and $d_{ij'}$ is equal to a $d_{i'j'}$ so can be flipped and absorbed.
In the first case, either $d_{ij'}\not\in J$, in which case it is  a $c_{ij'}$ and no
absorption is necessary, or $d_{ij'}\in J$ so that $d_{ij'}$ is
equal to a $d_{ij}$ and can be flipped and absorbed into 
$(u^*u)_{J- I'}$.

The proof for the third type of ``inner'' element, as well as the proofs 
 for the  ``outer'' elements are similar.

For the converse statement in (a), note that $c_{i'j}\in {I'}^c\cap J^c
\subset (I^c\cap J^c)\cup (I\cap J')\cup ({I'}^c\cap {J'}^c)$ 
and $d_{i'j}\in I'\cap J
\subset(I\cap J)\cup (I^c\cap {J'}^c)\cup (I'\cap J')$.
\qed

\medskip

With Lemma~\ref{lem:4.15}, the proof of (e) is completed. 

\medskip

{\bf Proof of (d):} 
If we let $w$ denote $u_{I,J}[u_{I',J'}]^*$, then
\[
w=[(uu^*)_{I- J}c_{ij}^1d_{ij}^1\cdots d_{ij}^sc_{ij}^{s+1}(u^*u)
_{J- I}]\ [(u^*u)_{J'- I'}c_{i'j'}^{r+1}d_{i'j'}^r
\cdots d_{i'j'}^1
c_{i'j'}^1(u^*u)_{I'- J'}].
\]

Since $I\ne I'$, there are two possibilities: either there exists $i_0
\in I- I'$ or there exists $i_0^\prime\in I'
- I$. We shall deal
with the first case only as the other is similar.

So assume first that
 $i_0\in I- I'$ and consider the two cases: $i_0\in J$
and $i_0\not\in J$.  In the first case $u_{i_0}$ is one of the
$d_{ij}$ and hence either $u_{i_0}$ is also a $c_{i'j'}$ in which case
$w=0$ by ``flipping'' and minimality; or $i_0\in J'- I'$ in which case
$w=0$ by ``hopping'' and minimality.

Now consider the case that $i_0\not\in J$.  In this case $i_0\in I- J$
and hence either $i_0\in J'$, in which case $i_0\in J'- I'$ and $w=0$
by ``hopping'' and minimality; or $i_0\not\in J'$, in which case $u_{i_0}$
is a $c_{i'j'}$ and $w=0$  again by ``hopping'' and minimality.

This proves the first statement in (d). The proof of the second statement
is achieved in a similar way.

The reader will note that ``maximality'' (meaning for instance that
$(u^*u)_J=0$ if $|J|>i_L$) was not used in the above proof
of (d). It's main use is in the proof of the important decomposition
(\ref{eq:sum}) below.

It being clear that (a) and (b) follow immediately from (d), it remains
to prove (c).

\medskip

{\bf Proof of (c):} 
In view of the strong orthogonality already proved, it will suffice to prove
that $u_{I,J}[u_{I,J}]^*u_{I',J}=u_{I',J}$ for $I\ne I'$.
We have 
\begin{eqnarray}\label{eq:colinear}\nonumber
u_{I,J}[u_{I,J}]^*u_{I',J}&=&[(uu^*)_{I- J}c_{ij}^1d_{ij}^1\cdots
d_{ij}^sc_{ij}^{s+1}(u^*u)_{J- I}]\\
&\times&
[c_{ij}^{s+1}d_{ij}^s\cdots
d_{ij}^1c_{ij}^1(uu^*)_{I- J}]\\\nonumber
&\times &
[(uu^*)_{I'- J}c_{i'j}^1d_{i'j}^1
\cdots d_{i'j}^rc_{i'j}^{r+1}(u^*u)_{J- I'}].
\end{eqnarray}

The term $(u^*u)_{J- I}$ in (\ref{eq:colinear}) can be absorbed
into the $d_{i'j}$'s or into $(u^*u)_{J- I'}$ by Lemma~\ref{lem:4.6}. Then in turn, the products
$c_{ij}^{s+1}(c_{ij}^{s+1})^*$, $d_{ij}^s(d_{ij}^s)^*,\ldots$, $c_{ij}^1
(c_{ij}^1)^*$ can be  alternatingly absorbed into  the combination of
$(uu^*)_{I'- J}$ and the $c_{i'j}$'s, or the combination of $(u^*u)_{J- I'}$ and the $d_{i'j}$'s.

 Finally both of the
occurences of the term $(uu^*)_{I- J}$ can also be absorbed
into either $(uu^*)_{I'- J}$ or a $c_{i'j}$,
 and what
remains is $u_{I',J}$. This completes the proof of 
Proposition~\ref{prop:4.14}.
\qed


\begin{defn}
In the special case of {\rm (\ref{eq:(I,J)})} where
 $I\cap J=\emptyset$, we have $s=0$ and $u_{I,J}$ 
has the form
\[
u_{I,J}=(uu^*)_Iu_c(u^*u)_J,
\]
where since $i_R+i_L=n+1$, 
 $I\cup J\cup \{c\}=\{1,\ldots,n\}$. 
We call such an element a {\rm ``one''}, and denote it by
$u_{I,c,J}$.
\end{defn} 

\begin{lem} For any $c\in \{1,\ldots,n\}$,
\begin{equation}\label{eq:sum}
u_c=\sum_{I,J}u_{I,J}=\sum_{I,J}u_{I,c,J}
\end{equation}
where the sum is taken over all disjoint $I,J$ satisfying {\rm (\ref{eq:12})}
and not containing $c$.
\end{lem}
\pf\
For convenience, let us say that for colinear partial isometries $u$ and
$v$, the formula $u=uv^*v+vv^*u$ is the result
of ``applying $v$ to $u$''. Given $c$, write $\{1,\ldots,n\}=\{c,
c_2,\ldots,c_n\}$.
The equation (\ref{eq:sum}) is obtained by
first applying $u_{c_2}$ to $u_c$, then applying $u_{c_3}$
to all occurrences of $u_c$, and in turn applying
$u_{c_4},\ldots,u_{c_n}$ to
all occurences of $u_c$ that are created in the previous step.

We thereby obtain
\[
u_c=\sum(uu^*)_Iu_c(u^*u)_J
\]
where the sum is over all disjoint subsets $I,J$ of 
$\{1,\ldots,n\}-\{c\}$ with 
\begin{equation}\label{eq:650}
I\cup J\cup\{c\}=\{1,\ldots,n\}.
\end{equation}
A term in this sum is zero unless $|I|\le i_R-1$ and $|J|\le i_L-1$. By
(\ref{eq:650}) and the fact that $i_R+i_L=n+1$,
$|I|=i_R-1$ and $|J|= i_L-1$ so (\ref{eq:sum}) follows. \qed

\medskip

Note that a change in 
the order of the $u_c$'s or $u_d$'s in (\ref{eq:(I,J)}) can at
most change the sign, since any such change can be accomplished by 
``flipping.'' In the next lemma, we consider elements defined by the right
side of (\ref{eq:(I,J)}) but without specifying an ordering of the $c$'s and
$d$'s. This lemma will enable us to define the signature $\epsilon(I,J)$ of
$u_{I,J}$ and prove the important Proposition~\ref{prop:4.17}.

\begin{lem}\label{lem:4.16}
Given $I,J$ with $|I|=i_R-1,|J|=i_L-1$, let $C=(I\cup J)^c$ and $D=I\cap J$.
For any permutations $(c_1,\ldots,c_{s+1})$ of $C$ 
and $(d_1,\dots,d_s)$ of $D$, the
element
\[
(uu^*)_{I- J}u_{c_1}u_{d_1}^*u_{c_2}u_{d_2}^*\cdots 
u_{c_s}u_{d_s}^*u_{c_{s+1}}
(u^*u)_{J- I}
\]
{\rm (}which equals $\pm u_{IJ}${\rm )}
decomposes uniquely as a product of ``ones'':
\begin{equation}\label{eq:19}
[u_{I_1,c_1,J_1}][u_{K_1,d_1,L_1}]^*
[u_{I_2,c_2,J_2}][u_{K_2,d_2,L_2}]^*\cdots 
[u_{I_s,c_s,J_s}][u_{K_s,d_s,L_s}]^*
[u_{I_{s+1},c_{s+1},J_{s+1}}],
\end{equation}
where the $I_i,J_j,K_k,L_l$ are uniquely determined by $I,J$ and
the $c$'s and $d$'s.
\end{lem}
\pf\
Let us first prove the existence. Each of the steps in the following equation
array
is achieved by ``expanding'' (for example, 
$u_{c_2}=u_{c_2}u_{c_2}^*u_{c_2}$) and/or
``hopping'':
\begin{eqnarray*}
\lefteqn{
(uu^*)_{I- J}u_{c_1}u_{d_1}^*u_{c_2}u_{d_2}^*
\cdots u_{c_s}u_{d_s}^*u_{c_{s+1}}(u^*u)_{J- I}}&\\
&=&(uu^*)_{(I- J)\cup (C- c_1)} u_{c_1}
[u_{d_1}^*u_{c_2}u_{d_2}^*\cdots
u_{c_s}u_{d_s}^*u_{c_{s+1}}]
(u^*u)_{(J- I)\cup (C- c_{s+1})}\\ 
&=&(uu^*)_{(I- J)\cup
(C- c_1)}u_{c_1}(u^*u)_{J- I} 
[u_{d_1}^*u_{c_2}u_{d_2}^*\cdots
u_{c_s}u_{d_s}^*u_{c_{s+1}}]
(u^*u)_{(J- I)\cup (C- c_{s+1})}\\
&=&
(uu^*)_{(I- J)\cup
(C- c_1)}u_{c_1}(u^*u)_{(J- I)\cup D} 
[u_{d_1}^*u_{c_2}u_{d_2}^*\cdots
u_{c_s}u_{d_s}^*u_{c_{s+1}}]
(u^*u)_{(J- I)\cup (C- c_{s+1})}\\
&=&(uu^*)_{(I- J)\cup (C- c_1)}u_{c_1}(u^*u)_{(J- I)
\cup D}\\
&\times&
[u_{d_1}^*u_{c_2}u_{d_2}^*
\cdots u_{c_s}u_{d_s}^*(uu^*)_{(I- J)\cup D}
u_{c_{s+1}}](u^*u)_{(J- I)\cup (C- c_{s+1})}\\
&=&[(uu^*)_{(I- J)\cup (C- c_1)}u_{c_1}
(u^*u)_{(J- I)\cup D}]\\
&\times&
[(u^*u)_ {(J- I)\cup\{c_1}\}u_{d_1}^*u_{c_2}u_{d_2}^*
\cdots u_{c_s}u_{d_s}^*(uu^*)_{(I- J)\cup\{c_{s+1}\}}]\\
&\times&
[(uu^*)_{(I- J)\cup D}
u_{c_{s+1}}(u^*u)_{(J- I)\cup (C- c_{s+1})}]
\end{eqnarray*}

This shows that $(uu^*)_{I- J}u_{c_1}u_{d_1}^*u_{c_2}
u_{d_2}^*\cdots u_{c_s}u_{d_s}^*u_{c_{s+1}}
(u^*u)_{J- I}$ equals
\[
[u_{(I- J)\cup (C- c_1),c_1,(J- I)\cup D}]
[u_{(I- J)\cup\{c_{s+1}\},(J- I)\cup\{c_1}\}]^*
[u_{(I- J)\cup D,c_{s+1},(J- I)\cup (C- c_{s+1})}],
\]
which is of the form $u_{I_1,c_1,J_1}[u_{I_2,J_2}]^*u_{I_3,c_{s+1},J_3}$,
so the existence follows by induction.

\medskip

We now prove the uniqueness. Look at the first three factors
of (\ref{eq:19}).
Since $u_{IJ}\ne 0$, by Proposition~\ref{prop:4.14} (d), we must have
$I_1=K_1$ and $L_1=J_2$. Furthermore, since $I_1\cup\{c_1\}\cup J_1=K_1\cup\{d_1\}\cup L_1$, we have
$J_2=(J_1\cup\{c_1\})-\{d_1\}$. Continuing, we see that all the 
sets $I_i,J_j,K_k,L_l$ are uniquely determined by  $J_1$ and
the $c$'s and $d$'s. A close look at  (\ref{eq:19}) and using
Proposition~\ref{prop:4.14} (d) and (e) reveals
that 
\[
u_{IJ}=\pm u_{I_1J_1}[u_{I_1J_{s+1}}]^*u_{I_{s+1}J_{s+1}}
\]
which equals $\pm u_{I_{s+1}J_1}$ by Proposition~\ref{prop:4.14} (e), 
so $u_{IJ}[u_{I_{s+1}J_1}]^*\ne 0$.
Then 
by Proposition~\ref{prop:4.14} (d) again,
$I=I_{s+1}$ and similarly $J=J_1$, completing the proof of uniqueness.
\qed


\begin{defn}
We assign a {\rm signature} to 
each ``one'' as follows: Let the  elements 
of $I$ be $i_1<i_2<\cdots<i_p$ {\rm (}where
$p=i_R-1${\rm )} and the elements of $J$ be $j_1<j_2<\cdots<j_q$ {\rm (}where 
$q=i_L-1${\rm )}. Then
$\epsilon(I,k,J)$ is defined to be the signature of
the permutation taking the $n$-tuple $(i_1,\ldots,i_p,k,j_1,\ldots,j_q)$ 
onto $(1,2,\ldots,n)$.

The {\rm signature} $\epsilon(I,J)$
 of $u_{I,J}$ is defined to be the product of the
signatures of the factors in its  decomposition {\rm (\ref{eq:19})} {\rm (}Recall 
that $u_{IJ}$ is defined so
that the $c$'s and $d$'s are in increasing order{\rm )}.
\end{defn} 

The next lemma will consider a 3-tuple
$(u_{I,J'},u_{I,J},u_{I',J})$, with $I\ne I'$ and $J\ne J'$, so that by
Proposition~\ref{prop:4.14}, $u_{I,J'}\perp u_{I',J}$,
$u_{I,J'}\top u_{I,J}$ and
$u_{I,J}\top u_{I',J}$.

Let us further assume that each
element of this 3-tuple is a ``one''. Then it is clear that 
$I'=(I-\{a\})\cup\{b\}$ and $J'=
(J-\{c\})\cup\{b\}$
for suitable elements $a\in I,\ c\in J$ and $b\in (I\cup J)^c$.
Hence
the 3-tuple has the form 
\begin{equation}\label{eq:19prime}
(u_{I,c,(J- c)\cup\{b\}}\ ,\ u_{I,b,J}\ ,\ 
u_{(I- a)\cup\{b\},a,J}).
\end{equation}

By direct calculation and simplification
\[
u_{IJ'}[u_{IJ}]^*u_{I'J}=(uu^*)_I\, u_c\, (u^*u)_J\, u_b^*\, 
(uu^*)_I\, u_a\, (u^*u)_J.
\]
Since $a\not\in J,b\not\in J$ and $I\cap J=\emptyset$, Lemma~\ref{lem:4.6}
shows that 
\[
u_{IJ'}[u_{IJ}]^*u_{I'J}=
(uu^*)_I\, u_cu_b^*\, (uu^*)_I\, u_a\, (u^*u)_J.
\]
Similarly, since $b\not\in I, c\not\in I$, we can remove the term $(uu^*)_I$ to obtain
\begin{equation}\label{eq:20}
u_{IJ'}[u_{IJ}]^*u_{I'J}=(uu^*)_I\, u_cu_b^*u_a\, (u^*u)_J,
\end{equation}
which equals $\pm u_{I'J'}$. 

Thus, from (\ref{eq:20}) and the uniqueness in Lemma~\ref{lem:4.16}, 
every such 3-tuple (\ref{eq:19prime}) of ``ones''
uniquely determines a 
corresponding 3-tuple of ``ones'' ($u_{I''J'},
u_{I''J''},u_{I'J''})$, such that
\begin{equation}\label{eq:24prime}
u_{I''J'}\, [u_{I''J''}]^*\, u_{I'J''}=\pm u_{I'J'}=
(uu^*)_I\, u_au_b^*u_c\, (u^*u)_J,
\end{equation}
where $I''=(I\cup\{c\})-\{a\}$ and $J''=(J\cup\{a\})-\{c\}$.
The given 3-tuple and the derived one  thus have the 
forms
\begin{equation}\label{eq:21}
(u_{I,c,(J- c)\cup\{b\}}\ ,\ u_{I,b,J}\ ,\ 
u_{(I- a)\cup\{b\},a,J})
\end{equation}
 and 
\begin{equation}\label{eq:21bis}
(u_{(I\cup\{c\})- a,a,(J- c)
\cup\{b\}}\ ,\ u_{(I\cup\{c\})- a,b,(J-
c)\cup\{a\}}\ ,\ u_{(I- a)\cup\{b\},c,(J- c)\cup\{a\}}).
\end{equation}

\begin{lem} \label{lem:6.4} Retain the above notation.
\begin{description}
\item[(a)] $u_{IJ'}[u_{IJ}]^*u_{I'J}=-u_{I''J'}[u_{I''J''}]^*u_{I'J''}$.
\item[(b)] $\epsilon(IJ')\epsilon(IJ)\epsilon(I'J)=-
\epsilon(I''J')\epsilon(I''J'')\epsilon(I'J'')$.
\item[(c)] For every 3-tuple {\rm (\ref{eq:19prime})} of ones,
\[
[\epsilon(IJ')u_{IJ'}]
[\epsilon(IJ)u_{IJ}]^*[\epsilon(I'J)u_{I'J}]=
[\epsilon(I''J')u_{I''J'}][\epsilon(I''J'')u_{I''J''}]^*
[\epsilon(I'J'')u_{I'J''}].
\]
\end{description}
\end{lem}
\pf\

(a) follows from (\ref{eq:20}) and (\ref{eq:24prime}).

(b) We shall use the more precise notation of (\ref{eq:21})
and (\ref{eq:21bis}).  Write
\[
u_{I,b,J}=(uu^*)_{\{i_1,i_2,\ldots,i_u,i_{u+1},\ldots,i_t,\ldots,i_{i_R-1}\}}
u_b^*(u^*u)_{\{j_1,j_2,\ldots,j_s,\ldots,j_v,j_{v+1},\ldots,j_{i_L-1}\}},
\]
where $i_u<c<i_{u+1}$, $i_t=a$, $j_s=c$ and $j_v<a<j_{v+1}$.

We can calculate $\epsilon(I,b,J)\epsilon(I'',b,J'')$ by counting the
number of transpositions required in ``moving'' $a$ from $I$ to $J''$
and $c$ from $J$ to $I''$. These are 
\begin{equation}\label{eq:11111}
(i_R-1)-t+1+(v+1)
\end{equation}
and
\begin{equation}\label{eq:22222}
(s-1)+1+(i_R-1)-u
\end{equation}
respectively.

We can calculate $\epsilon(I',a,J)\epsilon(I',c,J'')$ by counting the
number of transpositions required in ``moving'' $a$ from the ``middle''
 to $J''$
and $c$ from $J$ to the ``middle''. Taken together this is 
\begin{equation}\label{eq:33333}
v+(s-1).
\end{equation}

We can calculate $\epsilon(I,c,J')\epsilon(I'',a,J'')$ by counting the
number of transpositions required in ``moving'' $c$ from the ``middle''
 to $I''$
and $a$ from $I$ to the ``middle''. Taken together this is 
\begin{equation}\label{eq:44444}
(i_R-1)-u+(i_R-1-t).
\end{equation}

The sum of (\ref{eq:11111})-(\ref{eq:44444}) is $2(s+v-u-t)-1$ which is odd, so
exactly one or three of the numbers 
\[
\epsilon(I,b,J)\epsilon(I'',b,J'')\ ,\ 
\epsilon(I',a,J)\epsilon(I',c,J'')\ ,\ 
\epsilon(I,c,J')\epsilon(I'',a,J'')
\]
equals $-1$. In either case, (b) follows.

(c) follows from (a) and (b).
\qed

\begin{prop}\label{prop:4.17}
The family $\{\epsilon(IJ)u_{I,J}\}$ forms a rectangular grid which 
satisfies
\begin{equation}\label{eq:22}
\epsilon(IJ)u_{IJ}[\epsilon(IJ')u_{IJ'}]^*\epsilon(I'J')u_{I'J'}=
\epsilon(I'J)u_{I'J}.
\end{equation}
\end{prop}
\pf\
Since $u_{I'J'}[u_{IJ'}]^*=0$ for $I\ne I'$, 
the property (\ref{eq:703}) will follow from 
(\ref{eq:22}). The other grid properties are contained in
Proposition~\ref{prop:4.14}.

To prove (\ref{eq:22}), we apply Lemma~\ref{lem:4.16} to
decompose its left and right sides into ``ones.''
 To avoid cumbersome
notation let us denote any ``one'' with $u_c$ in the ``middle'' simply
by $(c)$ and its signature by $\epsilon$, with an identifying subscript.
With this convention, we have, for suitable $x_j,y_k,z_l,w_i\in
\{1,\ldots,n\}$ and $\epsilon_{pq}=\pm 1$,
\begin{eqnarray}\label{eq:sharp}\nonumber
\epsilon(IJ)u_{IJ}[\epsilon(IJ')u_{IJ'}]^*\epsilon(I'J')u_{I'J'}&=&
(\epsilon_{11} x_1)(\epsilon_{12} x_2)\cdots (\epsilon _{1,2r+1} x_{2r+1})\\
&\times&
[(\epsilon_{21} y_1)(\epsilon_{22} y_2)\cdots (\epsilon _{2,2s+1} x_{2s+1})
]^*\\\nonumber
&\times&
(\epsilon_{31} z_1)(\epsilon_{32} z_2)\cdots (\epsilon _{3,2t+1} z_{2t+1}),
\end{eqnarray}
and
\begin{equation}\label{eq:46prime}
\epsilon(I'J)u_{I'J}=
(\epsilon_{41} w_1)(\epsilon_{42} w_2)\cdots (\epsilon _{4,2n+1} w_{2n+1}).
\end{equation}

Recall that from Proposition~\ref{prop:4.14}(e), we have
\[
\epsilon(IJ)u_{IJ}[\epsilon(IJ')u_{IJ'}]^*\epsilon(I'J')u_{I'J'}=\pm 
u_{I'J}.
\]

Now each $w_i$ is either an $x_j,y_k$ or $z_l$ by the proof of 
Lemma~\ref{lem:4.15}, and
by that same Lemma, the $x_j,y_k,z_l$ which are not used,  call them
$v_1,\ldots,v_m$, occur twice with an even number of elements between them. 
Also, by Lemma~\ref{lem:6.4}(c) we may rearrange terms so that the right side of (\ref{eq:sharp}) becomes
\[
[(\epsilon_{51}v_1)(\epsilon_{51}v_1)
\cdots (\epsilon_{5m}v_m)(\epsilon_{5m}
v_{m})]\ 
[(\epsilon_{61} w_1)'(\epsilon_{62} w_2)'\cdots 
(\epsilon _{6,2n+1} w_{2n+1})'],
\]
where the notation $(\epsilon w)'$ indicates that the $I$ and
$J$ in $(w)$ may have changed. 

Since $(v_j)(v_j)^*(w_1)=(w_1)$ by colinearity, this collapses to
\begin{equation}\label{eq:46primeprime}
(\epsilon_{61} w_1)'(\epsilon_{62} w_2)'\cdots 
(\epsilon _{6,2n+1} w_{2n+1})'
\end{equation}

Since (\ref{eq:46prime}) and (\ref{eq:46primeprime}) are each equal to
$\pm u_{I'J}$, by the uniqueness in Lemma~\ref{lem:4.16}, 
$(\epsilon_{6i}w_i)'=(\epsilon_{4i}w_i)$, which proves (\ref{eq:22}). \qed

\section{Cartan factors of rank 1}

{\bf Proof of Proposition~\ref{prop:1} in the finite dimensional rank 1 case.}
It follows from Propositions 6.3 and 6.10 that the map $\epsilon(IJ)u_{IJ} \rightarrow E_{JI}$ is a ternary isomorphism onto
$n\choose i_R$ by $n\choose n-i_R+1$ complex matrices. By Remark~\ref{rem:6.05} and (\ref{eq:sum}),
$Y$ is completely 
isometric to a subtriple $H_n^{i_R}$ of a  Cartan factor
of type 1. In view of Lemma~\ref{lem:5.8} this completes the proof of Proposition~\ref{prop:1}  in the case that $Y$ is of
type 1 and rank 1 and finite dimensional. 

Note that the numbers $n$ and $k$ determine a simple
algorithm for constructing the unique matricial space $H_n^{k}$ (see the 
paragraph preceding Examples~1 and ~2 below).
The spaces $H_n^{k}$ are examples of the rank 1 JW*-triples whose
existence was assumed in  section 6, as the following lemma shows.
\begin{lem}
The spaces $H_n^{k}$ are rank 1 Hilbertian JC*-triples with $i_{R}=k$ and $i_{R}+i_{L}=n+1$. 
\end{lem}
\pf\
We will denote the generator $\sum_{I,J}\epsilon(IJ)E_{J,c,I}$ of the space $H_n^{k}$ by $u_{c}$. Note that the sum is
orthogonal by Proposition 6.3, so the $u_{c}$ are partial isometries. It is essential to notice that, for each $E_{J,c,I}$,
there are exactly $k$ (resp. $n-k+1$) elements $E_{J^{\prime},c^{\prime},
I^{\prime}}$ such that
$E_{J^{\prime},c^{\prime},I^{\prime}}[E_{J^{\prime},c^{\prime},I^{\prime}}]^*
E_{J,c,I}=E_{J,c,I}$ (resp.
$E_{J,c,I}[E_{J^{\prime},c^{\prime},I^{\prime}}]^*E_{J^{\prime},
c^{\prime},I^{\prime}}=E_{J,c,I}$), namely, those
$E_{J^{\prime},c^{\prime},I^{\prime}}$ with $(J^{\prime} - \{ c \})
 \cup \{c^{\prime}\} = J$ (resp. $(I^{\prime} -
\{c\}) \cup \{c^{\prime}\} = I$). In all other cases  
$[E_{J^{\prime},c^{\prime},I^{\prime}}]^*E_{J,c,I}=0$ (resp. $E_{J,c,I}
[E_{J^{\prime},c^{\prime},I^{\prime}}]^*=0$).
With this in mind, using Proposition 6.3, it is a straightforward verification to show that $\{u_{a} \,\ u_{a}
\,\ u_{b} \} = (1/2)u_{b}$, and $\{u_{a} \,\ u_{b} \,\ u_{a} \}=0$.  Lemma 6.9 (a) and the comments preceeding it together with
Proposition 6.3 shows easily that $\{u_{a} \,\ u_{b} \,\ u_{c} \}=0$. Hence, the $H_n^{k}$ are rank 1 JC*-triples and
are thus Hilbertian as discussed at the start of section 5.3. 

To see that $k=i_{R}$, consider the expression
\begin{equation}\label{ir}
u_{r}u_{r}^{\ast} \cdots u_{2}u_{2}^{\ast}u_{1}.
\end{equation}
If $r > k$, then, by the remarks above, for each term $\epsilon(IJ)E_{J,1,I}$ 
in the expansion of $u_{1}$, there must exist a
number $i$, $2 \leq i \leq r$ such that $u_{i}u_{i}^{\ast}\epsilon(IJ)E_{J,1,I}=0$. Hence, (\ref{ir}) is zero. Now assume $r=k$.
Suppose $I = \{2, \cdots, r\}$ and $J = \{r+1, \cdots, n\}$. Again by the above remarks, for each $2 \leq i \leq r$, there exists an
element $E_{J^{\prime},i,I^{\prime}}$ in the expansion of $u_{i}$ such that
$E_{J^{\prime},i,I^{\prime}}[E_{J^{\prime},i,I^{\prime}}]^*
E_{J,1,I}=E_{J,1,I}$, ensuring at least one nonzero term in
the expansion of (\ref{ir}). Since all possible nonzero terms of (\ref{ir}) are $\epsilon(IJ)E_{J,1,I}$ and those are
independent,
(\ref{ir}) is not zero. It follows that $i_{R}=k$. A similar argument shows that $i_{L}=n-k+1$.
\qed

\medskip

The following lemma implies the statement in Theorem 1 
that the $H_{n}^{k}$ are 1-mixed injectives.

\begin{lem}\label{ev} For each matrix $x=\sum a_{i}u_{i}$ in $H_{n}^{k}$,
\[
 tr( (xx^{\ast})^{1/2} )={n-1\choose k-1}^{1/2} (\sum
|a_{i}|^{2})^{1/2}.
\]
\end{lem}

\pf\
We first show that $xx^{\ast}$ can have at most one nonzero eigenvalue. Indeed, if $xx^*$ has at least two
distinct nonzero eigenvalues, then we may write
$xx^{\ast}=f(xx^{\ast})+g(xx^{\ast})$ for two nonzero disjointly supported even continuous functions $f$ and $g$ which vanish
at zero. Hence, $xx^{\ast}x=f(xx^{\ast})x+g(xx^{\ast})x$ is a non-trivial orthogonal decomposition of the element $xx^{\ast}x$
in the rank one JC*-triple $H_{n}^{k}$, which is impossible by definition of rank. Hence the eigenvalues of 
$xx^{\ast}$ are $\|x\|^{2}= \sum |a_{i}|^{2}$ and possibly zero. 

However, 
since each $u_{i}$ is the sum of exactly
$n-1\choose k-1$ orthogonal matrix units mutiplied by $\pm 1$, we have that
 $tr(xx^{\ast})={n-1\choose k-1} \sum |a_{i}|^{2}$  . Thus the multiplicity of the eigenvalue 
$\sum |a_{i}|^{2}$ is ${n-1\choose k-1}$ and $tr((xx^{\ast})^{1/2})={n-1\choose k-1}^{1/2} (\sum
|a_{i}|^{2})^{1/2}$.
\qed

\medskip

\begin{cor}The linear map $P$ defined by $Px=\sum tr(xu_{i}^{\ast}/{n-1\choose k-1}^{1/2})u_{i}$ is a contractive
projection from $n\choose k$ by $n\choose n-k+1$ complex matrices onto $H_{n}^{k}$.
\end{cor}
\pf\ 
Let $m$ denote the multiplicity $n-1\choose k-1$. Using Lemma \ref{ev} and the fact that the $H_{n}^{k}$ are Hilbertian, 
we see
that 

\begin{eqnarray*}
\|Px\|^{2} & = & \sum |tr(xu_{i}^{\ast}/ {m}^{1/2})|^{2}=tr(x(Px)^{\ast})/ {m}^{1/2}\\
 & \leq & \|x\|tr[(Px)(Px)^{\ast}]^{1/2}/ {m}^{1/2}\\
& =& \|x\| (\sum
|tr(xu_{i}^{\ast}/{m}^{1/2})|^{2})^{1/2}=\|x\|\|Px\| \qed
\end{eqnarray*}

\medskip

We now prove
 Proposition~\ref{prop:1} 
in the case that $Y$ is of type 1 and rank 1 and arbitrary dimension.
The key to the proof is the following lemma.

\begin{lem}\label{lem:719}
Suppose $Y$ is a \jwst\ of type 1 and rank 1 with grid $\{u_\lambda:\lambda
\in \Lambda\}$.  Then either $(uu^*)_I\ne 0$ for all finite subsets 
$I\subset\Lambda$, or $(u^*u)_J\ne 0$ for all finite subsets
$J\subset\Lambda$.
\end{lem}
\pf\
If $(u^*u)_I=0$ for some finite subset $I=\{i_1,\ldots,i_{n+1}\}$ we may
assume that
 $(u^*u)_{\{i_1,\ldots,i_n\}}\ne 0$.  Then as in (\ref{eq:719bis}),
\[
(u^*u)_{\{i_1,\ldots,i_n\}}=(u^*u)_{\{i_1,\ldots,i_n-1\}}
u_{i_n}^*(uu^*)_{\{i_n,i_n+1,\ldots,m\}}u_{i_n}
\]
for all $m\ge i_n$. Hence $(uu^*)_{\{i_n,i_n+1,\ldots,m\}}\ne 0$
for all $m\ge i_n$. Then by Lemma~\ref{lem:4.7}, $(uu^*)_J\ne 0$ for
all finite subsets $J\subset\Lambda$. \qed

\medskip

{\bf Proof of Proposition~\ref{prop:1} in the rank 1 type 1 case.}
We may assume $\dim(Y)=\infty$. For definiteness, we assume that
$(uu^{\ast})_{I} \ne 0$ for all finite 
subsets $I \subset \Lambda$. The other case in Lemma~\ref{lem:719}
is proved similarly. Let $E_{\lambda}$ denote $1 \otimes 
\psi_\lambda$ in
$B(H,\CC)$, where $\dim H=|\Lambda|$ and $\{ \psi_\lambda \}$ is an o.n.
 basis for $H$. By 
Proposition~\ref{prop:5.10bis}, for all finite subsets
$I\subset\Lambda$, the map
$\phi(u_{\lambda})=E_{\lambda}$ is a complete semi-isometry from $Y_I:=\mbox{sp}\, \{u_\lambda:\lambda
\in I\}$ to $sp\{ E_{\lambda}:\lambda \in I \}$. As a reflexive space, $Y$
is the norm-closure of the union of all the $Y_I$ as $I$ varies over
all finite subsets of $\Lambda$, so 
it follows that $Y$ is completely semi-isometric to 
$B(H,\CC)$. \qed

\medskip 

The proofs of Theorems~\ref{thm:1}, ~\ref{thm:2} and ~\ref{thm:3}(a)  being complete,
we now finish the proof of Theorem~\ref{thm:3}, give 
some examples of the spaces $H_n^k$, and pose some questions.

\medskip

{\bf Proof of Theorem~\ref{thm:3} for the rank 1 case.} 
Let $Y$ be an $n$-dimensional \jwst\ of rank 1.
It follows from Lemma~\ref{lem:5.8} and Proposition~\ref{prop:4.17}
that $Y=\mbox{Diag}\, (pY,(1-p)Y)$ where $pY$ and $(1-p)Y$
are triple isomorphic to $Y$, and $pY$ is completely isometric to 
some $H_n^{i_R}$. One now observes that the number $i_R$ for 
$(1-p)Y$ is strictly less than the $i_R$ for $Y$.  Indeed,
with $w_i=(1-p)u_i$, we have
\begin{eqnarray*}
(ww^*)_{\{1,2,\ldots,i_R\}}&=&(1-p)u_iu_i^*(1-p)u_2u_2^*(1-p)\cdots (1-p)
u_{i_R}u_{i_R}^*(1-p)\\
&=&(1-p)(uu^*)_{\{1,2,\ldots,i_R\}}\\
&=&\left(1-\sum_{|J|=i_R}(uu^*)_J\right)(uu^*)_{\{1,2,\ldots,i_R\}}\\
&=&(uu^*)_{\{1,2,\ldots,i_R\}}-(uu^*)_{\{1,2,\ldots,i_R\}}=0.
\end{eqnarray*}

Now set $Y_1=Y$, $p_1=p$, and $k_1=i_R$. Then setting
 $Y_2:=(1-p_1)Y_1$ and letting $k_2$
denote its $i_R$, then $k_2<k_1$.  Continuing in this way we see that
$Y=\mbox{Diag}\, (p_1Y_1,p_2Y_2,\ldots,p_mY_m)$, where each $p_jY_j$
is completely isometric to the space $H_n^{k_j}$. An application of
Lemma~\ref{lem:3.2} completes the proof of Theorem~\ref{thm:3}(b). \qed

\medskip

Note that the spaces
$\mbox{Diag}\, (H_n^{k_1},...,H_n^{k_m})$ are examples of Hilbertian rank 1 triples with
$i_R+i_L>n+1$, since $i_R=k_1$, $i_L=n-k_m+1$  and $k_1>...>k_m$.
Note also that the spaces $H_n^{i_R}$ can be explicitly constructed. Simply index columns (resp. rows) by combinations
$I$ (resp. $J)$ of
$\{ 1, \cdots ,n \}$ of length $i_R-1$ (resp. $i_L-1$).  Then define an orthonormal basis $\{ U_{i} \}$ for $H_n^{i_R}$ by the
requirement that $U_{i}$ equals the sum of all elements 
$\epsilon_{I,J} E_{J,I}$ where $I \cap J = \emptyset$ and $(I \cup J
)^{c}=i$ . Then choose signs $\epsilon_{I,J}$ by the procedure detailed above. We now give some
examples.

\begin{example}{\rm
Suppose that $Y=\mbox{sp}_{\CC}\, \{u_1,u_2,u_3\}$ and 
$i_R=i_L=2$. The rectangular grid given by
Proposition~\ref{prop:4.17} is depicted by the following array:
\[
\begin{array}{cc|ccc}
&&&I&\\
&&\{1\}&\{2\}&\{3\}\\\hline
&\{1\}&u_2u_1^*u_3&u_2u_2^*u_3u_1^*u_1&-u_3u_3^*u_2u_1^*u_1\\
J&\{2\}&-u_1u_1^*u_3u_2^*u_2&-u_1u_2^*u_3&u_3u_3^*u_1u_2^*u_2\\
&\{3\}&u_1u_1^*u_2u_3^*u_3&-u_2u_2^*u_1u_3^*u_3&u_1u_3^*u_2\\
\end{array}
\]

\medskip

By (\ref{eq:sum}) the ternary isomorphism from the span of this rectangular
grid to the canonical grid in $B(\CC^3)$, when restricted to $Y$, satisfies
\[
u_1\mapsto \left[\begin{array}{rrr}
0&0&0\\
0&0&1\\
0&-1&0
\end{array}\right]\quad ,\quad
u_2\mapsto \left[\begin{array}{rrr}
0&0&-1\\
0&0&0\\
1&0&0
\end{array}\right]\quad ,\quad
u_3\mapsto \left[\begin{array}{rrr}
0&1&0\\
-1&0&0\\
0&0&0
\end{array}\right].
\]

Thus $H_3^2$ is the subtriple of $B(\CC^3)$ consisting of all matrices of the
form
\[
\left[\begin{array}{rrr}
0&a&-b\\
-a&0&c\\
b&-c&0
\end{array}\right]
\]
and hence in this case $Y$ is actually completely semi-isometric to the 
Cartan factor $A(\CC^3)$ of 3 by 3 anti-symmetric complex matrices.
}
\end{example}

\begin{example} {\rm
Suppose that $Y=\mbox{sp}_{\CC}\, \{u_1,u_2,u_3,u_4\}$ and 
$i_R=3,i_L=2$. The rectangular grid given by
Proposition~\ref{prop:4.17} is depicted by the following array:
\[
\begin{array}{rr|rrrrrr}
&&&&I&&&\\ 
&&\{1,2\}&\{1,3\}&\{1,4\}&\{2,3\}&\{2,4\}&\{3,4\}\\ \hline
&\{1\}&22314&-33214&44213&-2233411&2244311&-3344211\\ 
J&\{2\}&11324&1133422&-1144322&33124&-44123&3344122\\ 
&\{3\}&-1122433&-11234&1144233&22134&-2244133&44132\\ 
&\{4\}&1122344&-1133244&11243&2233144&-22143&33142\\ 
\end{array}
\]

Here we have used the abbreviation $22314$ for $u_2u_2^*u_3u_1^*u_4$ and
so forth.

By (\ref{eq:sum}) the ternary isomorphism from the span of this rectangular
grid to the canonical grid in $B(\CC^6,\CC^4)$, when restricted to $Y$,
 satisfies
\[
u_1\mapsto \left[\begin{array}{rrrrrr}
0&0&0&0&0&0\\
0&0&0&0&0&1\\
0&0&0&0&-1&0\\
0&0&0&1&0&0
\end{array}\right]\quad ,\quad
u_2\mapsto 
\left[\begin{array}{rrrrrr}
0&0&0&0&0&-1\\
0&0&0&0&0&0\\
0&0&1&0&0&0\\
0&-1&0&0&0&0
\end{array}\right]\quad ,\quad
\]
and
\[
u_3\mapsto \left[\begin{array}{rrrrrr}
0&0&0&0&1&0\\
0&0&-1&0&0&0\\
0&0&0&0&0&0\\
1&0&0&0&0&0
\end{array}\right]\quad ,\quad
u_4\mapsto 
\left[\begin{array}{rrrrrr}
0&0&0&-1&0&0\\
0&1&0&0&0&0\\
-1&0&0&0&0&0\\
0&0&0&0&0&0
\end{array}\right]\quad ,
\]
so that  $Y$ is completely semi-isometric to $H_4^3$, which is the
subtriple of $B(\CC^6,\CC^4)$ consisting of all matrices of the form
\[
\left[\begin{array}{rrrrrr}
0&0&0&-d&c&-b\\
0&d&-c&0&0&a\\
-d&0&b&0&-a&0\\
c&-b&0&a&0&0
\end{array}\right].
\]
}
\end{example}

We now show that $H_3^2$ is not completely semi-isometric to $R_3$, as suggested to us by
N.\ Ozawa. It is clear that similar arguments can be used to prove Theorem~\ref{thm:1}(d).
Since $R_3$ is a homogeneous operator space, if there were a complete semi-isometry of $H^2_3$
onto $R_3$, then every isometry from $H^2_3$
onto $R_3$ would be a complete semi-isometry. In the notation of Example 1, let $U:H^2_3\rightarrow R_3\subset M_3(\CC)$ be the 
isometry defined by
\[
\left[\begin{array}{rrr}
0&0&0\\
0&0&1\\
0&-1&0
\end{array}\right]\mapsto
\left[\begin{array}{rrr}
0&0&1\\
0&0&0\\
0&0&0
\end{array}\right]
\quad ,\quad
\left[\begin{array}{rrr}
0&0&-1\\
0&0&0\\
1&0&0
\end{array}\right]\mapsto
\left[\begin{array}{rrr}
0&1&0\\
0&0&0\\
0&0&0
\end{array}\right],
\]
\[
\left[\begin{array}{rrr}
0&1&0\\
-1&0&0\\
0&0&0
\end{array}\right]\mapsto
\left[\begin{array}{rrr}
1&0&0\\
0&0&0\\
0&0&0
\end{array}\right].
\]
Then $U$ is not a complete contraction, since
\[
\left\|\left[
\begin{array}{rrrrrrrrr}
0&-1&0&0&0&-1&0&0&0\\
1&0&0&0&0&0&0&0&1\\
0&0&0&1&0&0&0&-1&0
\end{array}
\right]\right\|=\sqrt{2}
\]
and
\[
\left\|\left[
\begin{array}{rrrrrrrrr}
1&0&0&0&1&0&0&0&1\\
0&0&0&0&0&0&0&0&0\\
0&0&0&0&0&0&0&0&0
\end{array}
\right]\right\|=\sqrt{3}
\]


\begin{prob}
What is the
completely bounded Banach-Mazur distance $d_{cb}(H_n^k,R_n)$?
\end{prob}
\begin{prob}
What can one say about an arbitrary 1-mixed injective operator space?
What can one say about an arbitrary \jwst\ up to complete isometry?
\end{prob}

\begin{rem}
The authors hope to classify all 1-mixed injectives possessing a predual in a future publication by using the
known structure theory of $JBW^*$-triples in \cite{Horn} and \cite{HorNeh87}.
\end{rem}

\begin{rem}\label{rem:6.17}
After completing this paper, the authors discovered that the spaces $H_n^k$ appear, in a slightly different form,
 in \cite{AraFri78} in
their solution to the contractive projection problem on the
compact operators on a separable Hilbert space. Their
methods and proofs are different from ours.
In the special case that the projection is 
weak*-weak* continuous and $H$ is separable, Theorem 2 can be derived 
from their results. 
\end{rem}

\end{document}